\documentclass[11pt,a4paper,reqno]{amsart}
\usepackage{enumerate,array}
\usepackage[vcentermath,enableskew]{youngtab-op}
\usepackage{oppi-cb}
\newtheorem{theorem}{Theorem}[section]
\newtheorem{lemma}[theorem]{Lemma}
\newtheorem{corollary}[theorem]{Corollary}
\newtheorem{proposition}[theorem]{Proposition}
\newtheorem*{LRR}{\LR\ rule}
\newtheorem*{AC1}{Aldous's conjecture (v.1)}
\newtheorem*{AC2}{Aldous's conjecture (v.2)}
\newtheorem*{AC3}{Aldous's conjecture (v.3)}
\theoremstyle{definition}
\newtheorem{definition}[theorem]{Definition}

\theoremstyle{remark}
\newtheorem{remark}[theorem]{Remark}

\numberwithin{equation}{section}
\newcommand\yuung[1]{{\scriptsize$\young(#1)$}}
\newcommand\LR{{Littlewood--Richardson}}
\newcommand{\Adm}{\operatorname{Adm}}
\newcommand{\Adms}{\operatorname{Adm}^*}%
\newcommand{\ho}{\widehat\o}
\newcommand{\hg}{\widehat\g}
\newcommand{\hn}{\widehat\n}
\newcommand{\Tt}{t^\#}
\newcommand{\tshi}{\tilde t}
\newcommand\ubr[2]{\underset{#2}{\underbrace{#1}}}
\newcommand\leC{\le}
\newcommand\leCC{\ll}
\newcommand{\coimpdef}{\stackrel{\rm def}{\coimpp}}
\newcommand\lmax{{\l_{\rm max}}}
\newcommand\hh{\sigma}
\newcommand\eps{\varepsilon}
\newcommand{\bmdf}{{\ol B}}
\newcommand\repres[1]{\mathbf{#1}}
\newcommand\rT{\repres{T}}
\newcommand\rR{\repres{R}}
\newcommand\rY{\repres{Y}}
\newcommand\rD{\repres{D}}
\newcommand\secskip{\bigskip}
\begin{document} 
%

\title[Eigenvalues of Cayley graphs]{%
On the eigenvalues of Cayley graphs on the symmetric group
generated by a complete multipartite set of transpositions%
}
%
%
%
%
%
%
\author[Filippo Cesi]{Filippo Cesi}
\address{%
    Filippo Cesi\hfill\break
    \indent Dipartimento di Fisica\hfill\break
    \indent Universit\`a di Roma ``La Sapienza", Italy\hfil\break
    \indent and SMC, INFM-CNR.
}
\email{filippo.cesi@roma1.infn.it}

\subjclass[2000]{05C25, 05C50, 60K35}


%
%
%
\begin{abstract}
Given a finite simple graph $\cG$ with $n$ vertices,
we can construct the Cayley graph on the symmetric group $S_n$
generated by the edges of $\cG$, interpreted as transpositions.
We show that, if $\cG$ is complete multipartite, 
the eigenvalues of the Laplacian of $\Cay(\cG)$ have a
simple expression
in terms of the 
irreducible characters of transpositions,
and of the Littlewood-Richardson coefficients.
As a consequence we can prove that
the Laplacians of $\cG$ and of $\Cay(\cG)$ have
the same first nontrivial eigenvalue. 
This is equivalent to saying that Aldous's conjecture, asserting 
that the random walk and the interchange process have the
same spectral gap,
holds for complete multipartite graphs.
\end{abstract}

\maketitle
\thispagestyle{empty}
%
%
%
%
\vspace{5mm}
\section{Introduction} 

\noindent
Let $\cG = (V(\cG),E(\cG))$ be a finite graph
with $V(\cG)=\{1,2,\ldots,n\}$.
$\cG$ is always assumed to be simple, \iie\  
without multiple edges and loops, and undirected.
The \textit{Laplacian} of $\cG$ is the $n\times n$ matrix
$\D_\cG := D - A$, where $A$ is the adjacency matrix of $\cG$,
and $D = \diag(d_1,\ldots,d_n)$ with $d_i$ denoting
the degree of the vertex $i$.
Since $\D_\cG$ is symmetric and positive semidefinite,
its eigenvalues are real and nonnegative and can be ordered
as $0=\l_1 \le \l_2\le\cdots\le \l_n$. 
There is an extensive
literature dealing with bounds on the distribution
of the eigenvalues and consequences of these bounds.
We refer the reader to \cite{Big}, \cite{Chu} and \cite{GoRo}
for a general introduction to the subject. 
If (and only if) $\cG$ is connected, 
the second eigenvalue $\l_2$ is positive and is of particular
importance for several applications. 
The Laplacian $\D_\cG$ can be viewed as the generator%
\footnote{Or \textit{minus} the generator, depending 
on the preferred sign convention.}
of a continuous-time random walk on $V$, whose invariant measure
is the uniform measure $u$ on $V$.
In this respect,
$\l_2$ is the inverse of the ``relaxation time'' of the random
walk, a quantity related to the speed of convergence
to the invariant measure in $L^2(V,u)$ sense. $\l_2$ is also called the 
\textit{spectral gap} of $\D_\cG$.

There is a natural way to associate a Cayley graph to $\cG$.
Any edge $e=\{i,j\}$ in $E(\cG)$ (and, more generally,
any pair $\{i,j\}$ of elements of $V(\cG)$) can be identified
with a transposition $(ij)$ of the symmetric group $S_n$.
Consider then the Cayley graph with vertex set equal to $S_n$
and edges given by $(\pi,\pi e)$, where $\pi$ is a permutation
of $S_n$ and $e\in E(\cG)$. We let, for simplicity,
\begin{equation*}
  \Cay(\cG) := \Cay(S_n, E(\cG)) \,.
\end{equation*}
The Laplacian of $\Cay(\cG)$ is again the generator 
of a continuous-time Markov chain called the \textit{interchange
process} on $V$. It can be described as follows:
each site of $V$ is occupied by a particle of a different color, and
for each edge $\{i,j\}\in E(\cG)$, at rate 1, 
the particles at vertices $i$ and $j$ are exchanged.

It is easy to show (it follows from \eqref{eq:l2cay} and \eqref{eq:l2def},
but there are simpler and more direct proofs) that the spectrum of
$\D_\cG$ is a subset of the spectrum of $\D_{\Cay(\cG)}$.
By consequence
\begin{equation*}
  \l_2( \D_\cG) \ge \l_2(\D_{\Cay(\cG)}) \,.
\end{equation*}
Being an $n! \times n!$ matrix, in general the Laplacian of
$\Cay(\cG)$ has many more eigenvalues than the Laplacian of $\cG$.
Nevertheless, a neat conjecture due to David Aldous \cite{Ald}
states, equivalently:

\begin{AC1}
If $\cG$ is a finite connected simple graph, then   
\begin{equation*}
  \l_2( \D_\cG) = \l_2(\D_{\Cay(\cG)}) \,.
\end{equation*}
\end{AC1}

\begin{AC2}
If $\cG$ is a finite connected simple graph, then the random walk
and the interchange process on $\cG$ have the same spectral gap.
\end{AC2}

\smallno
Version 1 also appears in \cite{Fri2},
under the extra assumption of $\cG$ being bipartite.
Aldous's conjecture comes in a third flavor, which
originates from the analysis of the 
representations of the symmetric group.%
\footnote{This will be explained in greater detail in Sections
\ref{sec:rep} and \ref{sec:lap}.}
From this point of view the Laplacian $\D_\cG$ corresponds
to the $n$-dimensional \textit{defining} representation
of $S_n$, whose irreducible components are the trivial
representation and the representation associated with
the partition $(n-1,1)$. On the other hand
the Laplacian of the Cayley graph $\Cay(\cG)$ is associated
to the \textit{right regular} representation, which contains
\textit{all} irreducible representations. Let $\a=(\a_1, \ldots, \a_r)$
be a partition of $n$ and denote with $\rT^\a$ the irreducible
representation of $S_n$ which corresponds to the partition $\a$.
Let also $\lmax(\a)$ be the maximum eigenvalue of the matrix
\begin{equation*}
  \sum_{e=\{i,j\}\in E(\cG)} \rT^\a( e ) \,.
\end{equation*}
The trivial representation, which corresponds to the trivial
partition $\a = (n)$, is thus contained with multiplicity one
in both the defining and right regular representations, and
clearly $\lmax( (n) ) = |E(\cG)|$.
This accounts for the fact that the first eigenvalue
of both $\D_\cG$ and $\D_{\Cay(\cG)}$ is null.
By consequence, one finds (see Section \ref{sec:lap})
that Aldous's conjecture can be restated as follows:

\begin{AC3}
If $\cG$ is a finite connected simple graph, then   
\begin{equation*}
  \lmax(\a) \le \lmax((n-1,1)) \,,
\end{equation*}
for each nontrivial partition $\a$ of $n$, \iie\ 
for each partition $\a\ne (n)$.
\end{AC3}

\smallno
Aldous's conjecture has been proven for star-graphs in \cite{FOW}
and for complete graphs in \cite{DiSh}.
A major progress was made in \cite{HaJu} (similar results were
reobtained in \cite{KoNa}), where a rather
general technique was developed, which can be used to prove the
conjecture for trees (with weighted edges) and a few other cases.
Without entering into the details, we mention that this
technique is useful for classes of graphs whose spectral gap 
``tends to decrease'' when a new site, and relative edges, are
added to a preexisting graph. This is indeed the case of trees,
since it is not too difficult to prove that adding a leaf
with its relative edge cannot increase the spectral gap.
Using this approach, Aldous's conjecture has been recently
proven
\cite{Mor}, \cite{StCo} for hypercubes 
\textit{asymptotically}, \iie\  in the limit as the side length
of the cube tends to infinity.

The main result of the present paper is 
the proof of Aldous's conjecture
for \textit{complete multipartite graphs} (Theorem \ref{th:main}).
These are graphs such that
it is possible to write the vertex set of $\cG$ as a
disjoint union
\begin{equation*}
  V(\cG) = \{1,\, \ldots,\, n\} = N_1 \cup \cdots \cup N_p
\end{equation*}
in such a way that $\{i,j\}$ is an edge if and only if
$i$ and $j$ belong to distinct $N_k$'s.
The approach we follow, similar
in spirit to \cite{DiSh}, is group theoretical.

The plan of the paper is as follows.
After recalling a few standard facts on the representation
theory of the symmetric group in Section \ref{sec:rep},
we discuss the relationship between the Laplacian of Cayley graphs
and the irreducible representations of $S_n$ in Section \ref{sec:lap}.
In Section \ref{sec:out} the proof of our main result is outlined
in the case of bipartite graphs. Most of the relevant ideas
are discussed in this section.
Sections \ref{sec:pr} contains a detailed proof of the general
multipartite case. One of the key technical ingredients
is the identification of the \LR\ tableau with minimal
content, among all tableaux which appear in the decomposition
of a tensor product of representations of $S_n$
(Lemma \ref{th:lrmin}).
This aspect is discussed in Section \ref{sec:lr}.

After this paper was completed, a
beautiful proof of Aldous's conjecture has been found by Caputo, Liggett and
Richthammer \cite{CaLiRi}, which holds for \textit{arbitrary} graphs
(including weighted graphs).
Their approach is based on a subtle mapping, reminiscent of the
star-triangle transformation used in electric networks, which
allows a recursive proof.

\secskip
\section{The irreducible representations of $S_n$}
\label{sec:rep}

\noindent
We recall here a few well-known facts from the representation theory
of the symmetric group. The main purpose is to establish
our notation. Standard references for this section are
for instance \cite{CuRe} for general representation theory
and \cite{JaKe}, \cite{Sag} for the symmetric group.

Given a positive integer $n$, 
a \textit{composition} (resp. a \textit{weak composition}) of $n$
is a sequence $\a = (\a_1, \a_2, \a_3, \ldots)$
of positive (resp. nonnegative) integers  such that 
$\sum_{i=1}^\oo \a_i = n$. 
Since there is only a finite
number of nonzero terms, one can either consider
the whole infinite sequence or just the finite
sequence obtained by dropping all trailing zeros which
appear after the last nonzero element. 
We define the \textit{length} of $\a$ as the position 
of the last nonzero element in $\a$, so if the length of $\a$ is $r$,
we write
\begin{align*}
  \a = (\a_1, \ldots, \a_r) = (\a_1, \ldots, &\a_r, 0,0,0,\ldots)  \,.
\end{align*}
We also let $|\a| := \sum_{i} \a_i$, while, 
for an arbitrary set $S$, $|S|$ stands, as usual, for the
cardinality of $S$.
A \textit{partition} of $n$ is a nonincreasing composition of $n$.
We write $\a\scomp n$ if $\a$ is a composition of $n$,
$\a\wcomp n$ if $\a$ is a weak composition on $n$, and 
$\a\partit n$ if $\a$ is a partition of $n$.

We introduce a componentwise partial order 
in the set of all finite sequences of integers:
we write $\a \leC \b$
if $\a_i \le \b_i$ for each $i$.
If $\a,\b$ are partitions (compositions, weak compositions) we can define
the component-by-component sum $\a+\b$ which is still
a partition (composition, weak composition).
If $\a \leC \b$ the difference
$\b - \a$ is a weak composition of $|\b| - |\a|$.

\medno
The \textit{Young diagram} of a partition $\a$ of $n$ 
is a graphical representation of $\a$ as a collection
of $n$ boxes arranged in left-justified rows, with
the $i^{\rm th}$ row containing $\a_i$ boxes. For instance
\begin{equation*}
  (6,4,1) = \yng(6,4,1) \qquad\,.
\end{equation*}
We do not distinguish between a partition
and its associated Young diagram.
If the integer $k$ appears $m$ times in the partition $\a$
we may simply write $k^{m}$, so, for
instance,
\begin{equation*}
  (5,5,4,2,2,2,1,1) = (5^2, 4, 2^3, 1^2) \,.
\end{equation*}
Given a Young diagram $\a$, the \textit{conjugate}
(or \textit{transpose})
diagram is the diagram, denoted with $\a'$, obtained
from $\a$ by ``exchanging rows and columns'', for example
\begin{align*}
  \a  &= \yng(4,2,1) & 
  \a' &= \yng(3,2,1,1) \,. 
\end{align*}
The elements of $\a'$ are given by
\begin{equation}
  \label{eq:dual}
  \a'_s := |\{ j : \a_j \ge s \}| \,.
\end{equation}
Let $\Irr(S_n)$ be the set of all equivalence classes
of irreducible representations%
\footnote{In this paper, by ``representation'' we mean
a finite-dimensional representation over the field of complex numbers.}
of $S_n$. There is a one-to-one correspondence between
$\Irr(S_n)$ and the set of all partitions of $n$.
We denote with $[\a]$ the class of irreducible
representations of $S_n$ corresponding to the partition $\a$,
and with $f_\a$ the \textit{dimension} or \textit{degree} of
the representation.
For simplicity we write $[\a_1, \ldots, \a_r]$
instead of $[(\a_1, \ldots, \a_r)]$.
It is sometimes notationally convenient to refer to
a specific choice of a representative in the class $[\a]$.
We denote this choice with $\rT^\a$. Hence $\rT^\a$ is a
group homomorphism
\begin{equation*}
  \rT^\a : S_n \mapsto GL(f_\a,\bC) 
  \qquad \a\vdash n \,.
\end{equation*}
Since every
representation of a finite group is equivalent
to a unitary representation, we can assume (if useful) that $\rT^\a(\pi)$
is a unitary matrix for each $\pi \in S_n$.
Every representation $\rY$ of $S_n$
can be written, modulo equivalence,
as a direct sum of the $\rT^\a$
\begin{equation*}
  \rY \cong \bigoplus_{\a\partit n} \: c_\a \, \rT^\a \,.
\end{equation*}
A fundamental quantity associated to a representation $\rY$
of a finite group $G$ is the \textit{character} of $\rY$,
which we denote by $\c^{\rY}$, and is defined as
\begin{align*}
  \c^\rY &: G \mapsto \bC \\
       &: g \to \tr \rY(g) \,.
\end{align*}
Two representations are equivalent if and only if
they have the same characters, so, going back
to $G=S_n$, 
we can choose an arbitrary representative $\rT^\a$ in the
class $[\a]$ and define
$\c^\a := \c^{\rT^\a}$. The set $\{\c^\a : \a\partit n\}$
is the set of \textit{irreducible characters} of $S_n$.

If $H$ is a subgroup of $G$, we denote with $\rY\resR^G_H$
the \textit{restriction} of the representation $\rY$ to $H$.
Even if $\rY$ is an irreducible representation of $G$,
the restriction $\rY\resR^G_H$ is in general a reducible
representation of $H$. By consequence there exists
a collection of nonnegative integers $(c_{T})$ such that
\begin{equation*}
  \rY\resR^G_H = \bigoplus_{\rT \in \Irr(H)} c_{\,T} \, \rT \,.
\end{equation*}
When $G=S_n$ and $H$ is a Young subgroup 
$S_{(j,k)}$ (see Section \ref{sec:pr}) the coefficients
$c_{T}$ are called the \textit{\LR\ coefficients}.%
\footnote{The \LR\ coefficients are often 
equivalently (thanks to Frobenius reciprocity) defined
in terms of an \textit{induced} representation.}

\secskip
\section{Eigenvalues of Cayley graphs and representations of $S_n$} 
\label{sec:lap}

\noindent
We illustrate how, when studying the
eigenvalues of the Laplacian of Cayley graphs, one is (almost forcibly)
led to consider the irreducible representations
of the symmetric group. 
In this way we can show that version 3 of Aldous's conjecture
is equivalent to versions 1 and 2.
The material of this section
is more or less standard and overlaps with Section 4 of
\cite{DiSh}.

Given a finite set $S = \{s_1, \ldots, s_n\}$, we denote
with $\bC S$ 
the $n$-dimensional 
vector space which consists of all formal complex linear
combinations of the symbols $\{s_1\},\ldots,\{s_n\}$,
and with $\bC^S$ the vector space of
all functions $f:S \mapsto \bC$. 
$\bC^S$ is naturally isomorphic to $\bC S$ under the correspondence
\begin{equation}\label{eq:iso}
  f \otto \sum_{i=1}^n f(i) \, \{s_i\} \,.
\end{equation}
Any left action $(g,s)\to gs$ of a finite group $G$ on $S$
defines a representation%
$\rY$ of $G$ on $\bC S$ given by
\begin{equation}
\label{eq:Y}
  \rY(g) \Bigl( \sum_{i=1}^n a_i \, \{s_i\} \Bigr) :=
  \sum_{i=1}^n a_i \, \{ g s_i \} 
  \qquad g\in G\,.
\end{equation}
One can, equivalently, interpret $\rY$ as a representation on $\bC^S$,
in which case we have%
\footnote{With a slight abuse of notation we use the same 
symbol $\rY$ since the two representations 
are equivalent under \eqref{eq:iso}.}
\begin{equation}
\label{eq:Y1}
  [\rY(g) f](s) := f( g^{-1} s ) 
  \qquad g\in G,\ f\in\bC^S \,.
\end{equation}
$\bC G$ is the (complex) \textit{group algebra}
of $G$. Any representation $\rY$ of $G$ extends
to a representation of $\bC G$ by letting
\begin{equation*}
  \rY\Bigl( \sum_{g\in G} a_g \, g \Bigr) :=
  \sum_{g\in G} a_g \, \rY(g) 
  \qquad a_g\in\bC \,.
\end{equation*}
Let then $\cG$ be a finite graph with $V(\cG) = \{1,\ldots,n\}$.
The \textit{defining representation} of $S_n$, which
we denote by $\rD$, acts on $\bC V=\bC\{1,\ldots,n\}$ 
as
\begin{equation*}
  \rD(\pi)\Bigl( \sum_{i=1}^n a_i \, \{i\} \Bigr) =
  \sum_{i=1}^n a_i \, \{ \pi(i)\}
  \qquad \pi\in S_n \,.
\end{equation*}
The matrix elements of $\rD(\pi)$ in this basis are given by
\begin{equation*}
  [\rD(\pi)]_{ij} = 
  \begin{cases}
    1 & \text{if $j = \pi^{-1}(i)$} \\
    0 & \text{otherwise.}
  \end{cases}
\end{equation*}
The action of $\rD$ on the space $\bC^V$
is [$\rD(\pi) f](i) = f( \pi^{-1}(i))$, \iie\  $\rD(\pi) f = f \circ \pi^{-1}$.
If $\pi$ is a transposition, $\pi=(kl)$, we have
\begin{equation*}
[\rD\bigl((kl)\bigr)f](i) :=
\begin{cases}
  f(k) & \text{if $i=l$} \\
  f(l) & \text{if $i=k$} \\
  f(i) & \text{if $i\ne k,l$.}
\end{cases}
\end{equation*}
Hence, under the identification
of edges with transpositions of $S_n$
\begin{equation}
  \label{eq:ea}
  E(\cG) \ni e=\{i,j\} \too (ij) \in \{ \pi \in S_n : 
  \text{$\pi$ is a transposition} \}
\end{equation}
we can write
\begin{equation}\begin{split}
  \label{eq:lapg}
  &(\D_\cG f)(i) = \sum_{j:\, (ij)\in E(\cG) } \bigl[ f(i) - f(j) \bigr]
  \\
  &\quad
  = \sum_{j:\, e=(ij)\in E(\cG) } \bigl[ f(i) - \rD(e) f(i) \bigr] 
  =
  \sum_{e \in E(\cG) } \bigl[ f(i) - \rD(e) f(i) \bigr] \,,
\end{split}\end{equation}
where, in the last term, we have included
the null contribution of those edges with both
endpoints different from $i$. The reason is that we can now
rewrite \eqref{eq:lapg} in operator form.
If denote with $\id_n$ the identity operator
acting on an $n$-dimensional vector space, we have
\begin{equation}
  \label{eq:lap=}
  \D_\cG = |E(\cG)| \, \id_n - \sum_{e\in E(\cG)} \rD(e) 
  = |E(\cG)| \, \id_n - \rD \Bigl( \sum_{e\in E(\cG)} e \Bigr) \,,
\end{equation}
where, in view of correspondence \eqref{eq:ea}, $\sum_{e\in E(\cG)} e$ can
be considered an element of the group algebra $\bC S_n$,
and, in the last equality, we have used the linear
extension of $\rD$ to a representation of $\bC S_n$. 
Given a finite graph $\cG$ we define
\begin{equation}
  \label{eq:WG}
  W(\cG) := \sum_{e\in E(\cG)} e \in \bC S_n
\end{equation}
and rewrite \eqref{eq:lap=} as 
\begin{equation}
  \label{eq:lap=2}
  \D_\cG = |E(\cG)| \, \id_n - \rD[ W(\cG)] \,.
\end{equation}
A relationship for the corresponding eigenvalues
trivially follows
\begin{equation}
  \label{eq:lap=e}
  \l_i(\D_\cG) = |E(\cG)| - \l_{n-i}\bigl( \rD[W(\cG)] \bigr)
  \qquad i =1,\ldots,n \,. 
\end{equation}
We remark that, in the more general case of a \textit{weighted}
graph with edge weights $(w_e)_{e\in E(\cG)}$, identities
\eqref{eq:lap=2} and \eqref{eq:lap=e} remain valid as long as
one uses the ``correct definition'' of 
$W(\cG)$ as $W(\cG):= \sum_{e\in E(\cG)} w_e \, e$.

\smallno
We can associate to the graph $\cG$ the Cayley graph
$\Cay(S_n, E(\cG))$ with vertex set $S_n$, where $n$ is the cardinality
of $V(\cG)$, and edge set given by 
\begin{equation*}
  \{ (\pi, \pi e) : \pi\in S_n,\ e=(ij) \in E(\cG) \} \,.
\end{equation*}
Since each transposition coincides with its inverse, 
this Cayley graph is undirected.
We let for simplicity $\Cay(\cG) := \Cay(S_n, E(\cG))$.
If we denote with $\rR$ the \textit{right regular representation}
of $S_n$ which acts on $S_n$ and on $\bC^{S_n}$ respectively as%
\footnote{The right regular representation is a \textit{left}
action, like every representation.}
\begin{align*}
  \rR(\pi)  \pi' &= 
  \pi' \, \pi^{-1} 
  &
  &\qquad \pi, \pi'\in S_n, \\
  [\rR(\pi) f](\pi') &= f( \pi' \pi)
  &
  &\qquad f: S_n\mapsto \bC \,,
\end{align*}
we can proceed as in \eqref{eq:lapg} and obtain
\begin{equation}\begin{split}
  \label{eq:lapcay}
  (\D_{\Cay(\cG)} f)(\pi) &= \sum_{e \in E(\cG) } 
  \bigl[ f(\pi) - f(\pi e) \bigr]
  = \sum_{e \in E(\cG) } \bigl[ f(\pi) - \rR(e) f(\pi) \bigr] \,.
\end{split}\end{equation}
Identities \eqref{eq:lap=2} and \eqref{eq:lap=e} become, for the Cayley graph,
\begin{align}
  \label{eq:lc}
  \D_{\Cay(\cG)} &= |E(\cG)| \, \id_{n!} - \rR[ W(\cG)] \\
  \label{eq:lce}
  \l_i\bigl( \D_{\Cay(\cG)} \bigr) 
    &= |E(\cG)| - \l_{n!-i}\bigl( \rR[W(\cG)] \bigr)
  \qquad i =1,\ldots,n! \,. 
\end{align}
The right regular representation $\rR$ is equivalent to the
left regular representation (under the change of basis 
$\pi \to \pi^{-1}$)
and can be written as
a direct sum of \textit{all} irreducible representations,
each appearing with a multiplicity equal to its dimension
\begin{equation*}
  [\rR] = \bigoplus_{\a \partit n} f_\a\, [\a] \,.
\end{equation*}
By consequence
the spectrum of $\rR[W(\cG)]$ can be written as%
\footnote{To get the correct multiplicities of the eigenvalues
one must include the coefficients $f_\a$ and interpret
the union over $\a$ as a disjoint union of multisets.}
\begin{equation}
  \label{eq:spec}
  \spec \rR[W(\cG)] = \bigcup_{\a\partit n} \spec \rT^\a[ W(\cG) ] \,.
\end{equation}
The (trivial) one-dimensional identity representation 
$\rT^{(n)} = I_1$, corresponding
to the partition $(n)$, appears in this decomposition
exactly once and we have 
$\rT^{(n)} [W(\cG)] = |E(\cG)| \cdot I_1$; thus
its unique eigenvalue is equal to $|E(\cG)|$, which accounts for
the fact that $\l_1\bigl( \D_{\Cay(\cG)} \bigr) = 0$.
If $\cG$ is connected, the set $E(\cG)$, considered as
a set of transpositions,  generates $S_n$, and hence
$\Cay(\cG)$ is also connected and $\l_1$ is the
unique null eigenvalue of $\D_{\Cay(\cG)}$. In \textit{any} case,
letting
\begin{equation*}
  \lmax\bigl[ \a,\, W(\cG) \bigr] := \max \spec \rT^\a[ W(\cG) ] =
  \l_{f_\a}\bigl( \rT^\a[ W(\cG) ] \bigr)\,,
\end{equation*}
we have, for what concerns the second eigenvalue of the Cayley
graph,
\begin{equation}
  \label{eq:l2cay}
  \l_2\bigl( \D_{\Cay(\cG)} \bigr) = 
  |E(\cG)| - \max_{\a\partit n,\, \a\ne(n)} \lmax\bigl[\a,\, W(\cG)\bigr] \,.
\end{equation}
On the other hand, the defining representation can be decomposed as
$[\rD] \penalty0 = [n] \oplus \penalty0 [n-1,1]$,
which implies
\begin{equation}
  \label{eq:l2def}
  \l_2\bigl( \D_\cG \bigr) = 
  |E(\cG)| - \lmax\bigl[ (n-1,1),\, W(\cG) \bigr] \,.
\end{equation}
The main result of paper is the following:

\begin{theorem}\label{th:main}
If $\cG$ is a complete multipartite graph with $n$ vertices, we have
\begin{equation}\label{eq:main}
  \lmax\bigl[ \a,\, W(\cG) \bigr]  \le
  \lmax\bigl[ (n-1,1),\, W(\cG) \bigr] 
\end{equation}
for all irreducible representations $[\a]$ of $S_n$ with $[\a] \ne [n]$.
\end{theorem}

\smallno
From \eqref{eq:l2cay}, \eqref{eq:l2def}, and Theorem \ref{th:main}
it follows that:

\begin{corollary}\label{th:equ}
If $\cG$ is a complete multipartite graph, then Aldous's conjecture
holds, that is
\begin{equation*}
  \l_2\bigl( \D_{\Cay(\cG)} \bigr) =   \l_2\bigl( \D_\cG \bigr) \,.
\end{equation*}
\end{corollary}

\secskip
\section{Outline of the proof in the bipartite case} 
\label{sec:out}

\noindent
We briefly sketch in this section the proof of Theorem \ref{th:main}
in the bipartite case, which requires less notation
than the more general multipartite case but illustrates
most of the relevant ideas. The general case can be
treated by relatively standard induction. All missing details
will be found in later sections.

\medno
We start with a well-known fact \cite{DiSh} 
about the complete graph $K_n$.
Given an irreducible representation $[\a]$ of $S_n$,
corresponding to the partition $\a = (\a_1, \ldots, \a_r)$,
we consider the normalized character on the sum of
all transpositions
\begin{equation*}
  q_\a := \frac{n(n-1)}{2 f_\a} \, \c^\a( e ) 
  \qquad \a\partit n\,,
\end{equation*}
where $e$ is an arbitrary transposition of $S_n$. 
In the case of transpositions,
Frobenius formulas for the irreducible characters
take the simple form \cite{Ing}
\begin{equation}
  \label{eq:q=}
  q_\a = \ov2 \sum_{i=1}^\oo \a_i \, [ \a_i - (2i-1)] 
  = \ov2 \sum_{i=1}^r \a_i \, [ \a_i - (2i-1)] \,,
\end{equation}
where $r$ is the length of $\a$.
We use expression \eqref{eq:q=} as a \textit{definition} of $q_\a$
when $\a$ is, more generally,  
a \textit{weak composition} of $n$, even though, when
$\a$ is not a partition, the quantity $q_\a$ has no significance
associated to an irreducible representation of $S_n$.
A simple application of the Schur's lemma yields
the following result
(see \cite[Lemma 5]{DiSh} for a more
general statement where arbitrary conjugacy classes are considered).

\begin{proposition}\label{th:kn}
If $\,\rT^\a$ is an irreducible representation of $S_n$
corresponding to the partition $\a\partit n$, then
\begin{align*}
  \rT^\a[W(K_n)] = q_\a \, \id_{f_\a} \,.
\end{align*}
\end{proposition}

\medno
Let then $n=j+k$ with $j,k$ two positive integers
and consider the complete bipartite graph $K_{j,k}$ with 
vertex set $\{1,\ldots,n\}$
and edges 
$\{i,i'\}$ with $i \le j$ and $i'>j$.
Since the complement of $K_{j,k}$ is given by
\begin{equation*}
  \ol K_{j,k} = K_j \cup K_k \,,
\end{equation*}
using Proposition \ref{th:kn}, one can prove (see Proposition \ref{th:keta})
that the eigenvalues of $\rT^\a[W(K_{j,k})]$
have the form
\begin{equation}
  \label{eq:form}
  q_\a - q_\b -q_\g \,,
\end{equation}
where $\b$ is a partition of $j$ and $\g$ is a partition of
$k$, subject to the condition that the \LR\ coefficient $c^\a_{\b,\g}$
is positive. The reason for this is that
the irreducible representation $[\a]$ of $S_n$
is no longer irreducible when restricted to the Young
subgroup $S_{(j,k)}\cong S_j \times S_k$, but
it is a direct sum of irreducible components
\begin{equation}
  [\a] \resR^{S_n}_{S_{(j,k)}} = 
  \bigoplus_{\b\partit j, \, \g \partit k} c^\a_{\b,\g} \; [\b]\otimes[\g] \,.
  \label{eq:lr}
\end{equation}
If one is interested in keeping track of multiplicities, each
pair $(\b, \g)$ appearing in \eqref{eq:lr} contributes with a multiplicity
equal to $c^\a_{\b,\g} \: f_\b f_\g$.
For example, using the \LR\ 
rule (see Section \ref{sec:lr}), we find
the decomposition
\begin{align*}
  &[4,2,1] \resR^{S_7}_{S_{(4,3)}} =
  [4] \otimes [2,1] \oplus
  [3,1] \otimes [3] \oplus
  2\; [3,1] \otimes [2,1] 
  \oplus
  [3,1] \otimes [1^3] 
  \\
  &\quad
  \oplus
  [2^2] \otimes [3] \oplus
  [2^2] \otimes [2,1] \oplus
  [2,1^2] \otimes [3] \oplus
  [2,1^2] \otimes [2,1] \,.
\end{align*}
This, in turn, determines that the eigenvalues of $\rT^{(4,2,1)}[W(K_{4,3})]$
are those given in Table \ref{tab:exa}. 
Thus $\lmax( (4,2,1), \, W[K_{4,3}] ) = 5$.

\begin{table}
\begin{center}
\caption{Eigenvalues of $\rT^{(4,2,1)}[W(K_{4,3})]$}
\label{tab:exa}%
\renewcommand\arraystretch{1.2}
\begin{tabular}{||c|c|c||c|c|c||} \hline 
\multicolumn{6}{|c|}{%
$\vphantom{\Bigl[}%
\a=(4,2,1) \quad \l = q_\a - q_\b - q_\g$} 
\\ \hline
\vphantom{\Large g} $\b$   &  $\g$   &   $\l$  & $\b$   &  $\g$   &   $\l$
\\ \hline
$(4)$  & $(2,1)$ & $-3$   &
$(3,1)$  & $(3)$ & $-2$   \\ \hline
$(3,1)$  & $(2,1)$ & $1$   &
$(3,1)$  & $(1^3)$ & $4$   \\ \hline
$(2,2)$  & $(3)$ & $0$   &
$(2,2)$  & $(2,1)$ & $3$   \\ \hline
$(2,1,1)$  & $(3)$ & $2$   &
$(2,1,1)$  & $(2,1)$ & $5$   \\ \hline
\end{tabular}
\end{center}
\end{table}

We say that the pair $(\b,\g)$ is $\a$--admissible if
$c^\a_{\b,\g}>0$ and we define
\begin{equation*}
  B^\a_{j,k} := 
  \max_{\text{$\a$--admissible $(\b,\g)$}} 
  q_\a - q_\b -q_\g \,,
\end{equation*}
so that
\begin{equation*}
    \lmax\bigl[\a, W(K_{j,k})\bigr] = B^\a_{j,k} \,.
\end{equation*}
In general, given $\a\partit n$ and $\b\partit j$,
there are several different $\g\partit k$ such that
$(\b,\g)$ is $\a$--admissible. One of the central points
of the proof is the identification of the particular
$\hg=\hg(\a,\b)$ which corresponds to a minimum value of $q_\g$,
given $\a$ and $\b$, so that
\begin{equation}
  \label{eq:Bq}
  B^\a_{j,k} = 
  \max_{\b\partit j,\; \b\le \a } 
  q_\a - q_\b -q_{\hg(\a,\b)} 
  \,,
\end{equation}
What we find, in particular, is that (Lemma \ref{th:lrmin})
\begin{equation*}
  \hg(\a,\b) = \srt(\a-\b)
\end{equation*}
where $\srt$ is the operator that sorts a sequence
in nondecreasing order in such a way
that the resulting sequence is a partition.
So, for instance, if $\a=(7,6,2,1)$ and $\b=(5,2,2)$
we have
\begin{equation*}
  \a-\b = (2,4,0,1) \qquad
  \hg= \srt(\a-\b) = (4,2,1) \,.
\end{equation*}
At this point one could reasonably hope in  some monotonicity
property of the $B^\a_{j,k}$ with respect to $\a$.
There is a partial order ``$\leT$'' in the set of all
partitions of $n$, called \textit{dominance}
(see Section \ref{sec:pr}), which plays a crucial
role in the representation theory of the symmetric group.
It would be nice to prove something like
\begin{equation}
  \label{eq:blb}
  \a \leT \a' 
    \quad \impp \quad
    \lmax\bigl[\a,\, W(K_{j,k})\bigr] \le
    \lmax\bigl[\a',\, W(K_{j,k})\bigr] 
  \,.
\end{equation}
Since any nontrivial partition $\a$
of $n$ is dominated by the partition $(n-1,1)$, 
property \eqref{eq:blb}, if true, 
would imply Theorem \ref{th:main} for $\cG=K_{j,k}$.
Implication
\eqref{eq:blb} is unfortunately false.%
\footnote{A simple counterexample is given by $K_{3,1}$.
One easily finds that, in this case, we have
$\lmax\bigl[ (2,1,1),\, W(K_{3,1}) \bigr] = 1$ and
$\lmax\bigl[ (2,2),\, W(K_{3,1}) \bigr] = 0$.}
Nevertheless,
our actual strategy is a slight detour from this monotonicity idea.
We consider a modified version of the quantities \eqref{eq:Bq}
\begin{equation}
  \label{eq:Bq1}
  \bmdf^\a_{j,k} := 
  \max_{\b\partit j,\; \b\le \a } 
  q_\a - q_\b -q_{\a-\b} 
  \,.
\end{equation}
Then we realize (Proposition \ref{th:qinc})
that $q_{\a-\b} \le q_{\srt(\a-\b)}$, and thus, by consequence,
$ B^\a_{j,k} \le  \bmdf^\a_{j,k}$.
Using \eqref{eq:q=} one finds (Proposition \ref{th:easy})
a very simple expression for the quantity $q_\a -q_\b - q_{\a-\b}$, namely
\begin{equation*}
  q_\a -q_\b - q_{\a-\b} = \b \cdot (\a-\b) = 
  \sum_{i=1}^\oo \b_i \: (\a_i-\b_i) \,.
\end{equation*}
At this point one gets a lucky break. In fact
\begin{enumerate}[(a)]
\item 
The monotonicity property \eqref{eq:blb} holds for the quantities
$\bmdf^\a_{j,k}$ (Proposition \ref{th:mono}).
\item
If $\a=(n-1,1)$, we find that%
\footnote{Unless $j=k=1$, but this case is trivial.}
$\bmdf^{(n-1,1)}_{j,k} = B^{(n-1,1)}_{j,k}$%
(Proposition \ref{th:qqq}).
\end{enumerate}
Combining these facts, we obtain, for any $\a\partit n$ with $\a\ne(n)$,
\begin{equation*}
  B^\a_{j,k} \le \bmdf^\a_{j,k} 
  \le \bmdf^{(n-1,1)}_{j,k} = B^{(n-1,1)}_{j,k} \,,
\end{equation*}
and Theorem \ref{th:main} is proven.

\secskip
\section{Proof of Theorem \ref{th:main}} 
\label{sec:pr}

\noindent
A complete multipartite graph with $n$ vertices is identified,
up to a graph isomorphism, by a partition of $n$, so,
if $\hh = (\hh_1, \ldots, \hh_p)\vdash n$, with $p\ge 2$,
we denote
the associated complete multipartite graph 
with $K_\hh = K_{\hh_1,\ldots,\hh_p}$.
The set $\{1,\ldots,n\}$ can be written as a disjoint union
\begin{equation*}
  \{1, \ldots, n\} = N^\hh_1 \cup \cdots \cup N^\hh_p 
\end{equation*}
of subsets $N^\hh_k$ of cardinality $\hh_k$ given by
\begin{equation}
  \label{eq:N}
  N^\hh_k: = \{\hh_1+\cdots+\hh_{k-1}+1\, ,
  \; \ldots \;,\, \hh_1+\cdots+\hh_k\} \,.
\end{equation}
Let $S^\hh_{k}$ be the subgroup of $S_n$ which consists of the
permutations $\pi$ such that $\pi(i) = i$ for each 
$i\in \{1,\ldots,n\}\setm N^\hh_k$.
The \textit{Young subgroup} $S_\hh$ is defined as
\begin{equation*}
  S_\hh = S_{(\hh_1,\ldots,\hh_p)} := 
  S^\hh_1 \times \cdots \times S^\hh_p \,.
\end{equation*}
In other word a permutation $\pi$ belongs to $S_\hh$ if and only if
\begin{equation}
  \label{eq:ysub}
  i\in N^\hh_k \quad\impp\quad \pi(i) \in N^\hh_k \,.
\end{equation}
The subgroup $S_\hh$ is naturally isomorphic to the
(exterior) Cartesian product
$S_{\hh_1} \times \cdots\times S_{\hh_p}$. 

\smallno
We observe that the complement of $K_\hh$ is a disjoint union
of complete graphs
\begin{equation*}
  \ol K_\hh = K_{\hh_1} \cup \cdots \cup K_{\hh_p} \,,
\end{equation*}
and hence
\begin{equation*}
  W[K_\hh]  = W[K_n] - W[ \cup_{k=1}^p K_{\hh_k})] 
\end{equation*}
and thanks to Proposition \ref{th:kn}, we get, for
any irreducible representation $[\a]$ of $S_n$, the
identity%
\footnote{
Even though $\a$ and $\hh$ are both partitions of $n$,
they play a very different role. $[\a]$ is an 
equivalence class of irreducible
representations of $S_n$, while $\hh$ determines the
structure of the graph $K_\hh$.}
\begin{align}
  \label{eq:tq}
  \rT^\a [W(K_\hh)] = q_\a \id_{f_\a} - \rT^\a[ W( \cup_{k=1}^p K_{\hh_k})] \,.
\end{align}
The quantity $W( \cup_{k=1}^p K_{\hh_k})$ belongs to
the group algebra of the Young subgroup $S_\hh$.
The irreducible representation $[\a]$ of $S_n$
is no longer irreducible when restricted to $S_{\hh}$.
The irreducible representations of 
$S_{\hh} \cong S_{\hh_1} \times \cdots\times S_{\hh_p}$
are in fact the (outer) tensor products of the irreducible
representations of each $S_{\hh_i}$%
\begin{equation*}
  \Irr(S_\hh) = \bigl\{\, [\b^1] \otimes\cdots\otimes [\b^p] \;:\:
  \b^i \partit \hh_i \text{ for each $i=1,\dots,p$} \bigr\}\,.
\end{equation*} 
The obvious step at this point is to take advantage of the decomposition
\begin{equation}
  [\a] \resR^{S_n}_{S_{\hh}} = 
  \bigoplus_{\b^1 \partit \hh_1,\, \ldots,\, \b^p\partit \hh_p} 
  c^\a_{\b^1,\ldots,\b^p} \; [\b^1]\otimes\cdots\otimes [\b^p] 
  \label{eq:lrp}
\end{equation}
into a sum of irreducible representations of $S_\hh$.
Identities \eqref{eq:tq}, \eqref{eq:lrp} and the fact that
$W( \cup_{k=1}^p K_{\hh_k}) \in \bC S_\hh$,
imply that the eigenvalues
of $\rT^\a[W(\cG)]$ are of the form $q_\a - \l$, where
$\l$ is an eigenvalue of
\begin{equation}
  \label{eq:tbe}
  \bigotimes_{k=1}^p \rT^{\b_k} [ W( \cup_{k=1}^p K_{\hh_k}) ] \,,
\end{equation}
and where $(\b^1, \ldots, \b^p)$ 
is a collection of partitions $\b^k\partit \hh_k$
such that the (multi) \LR\ coefficient $c^\a_{\b^1,\ldots,\b^p}$
is positive. We give a name to these collections of $\b^k$.

\begin{definition}\label{th:adm}
Given the partitions $\a \vdash n$ and 
$\hh=(\hh_1,\ldots,\hh_p) \vdash n$, we say that the $p$-tuple
$(\b^1, \ldots, \b^p)$ of partitions is
\textit{$(\a,\hh)$--admissible} if
\begin{enumerate}[(i)]
\item each $\b^i$ is a partition of $\hh_i$
\item $c^\a_{\b^1,\ldots,\b^p} > 0$.
\end{enumerate}
We denote with $\Adm(\a,\hh)$ the set of all $(\a,\hh)$--admissible
$p$-tuples of partitions.
\end{definition}

\medno
The spectrum of the matrix in \eqref{eq:tbe}
can be expressed in a simple form thanks to the fact
that $\cup_{i=1}^p K_{\hh_i}$ is a disjoint union.

\begin{proposition}\label{th:prod}
Let $\cH$ be a finite graph which is the (disjoint) union of
subgraphs $\cH_1, \ldots, \cH_p$, and let $\hh_i$ be the number
of vertices of $\cH_i$.
For each $i=1,\ldots,p$, let $\rY^i$ be a representation
of $S_{\hh_i}$ and let $\rY$ be the representation of $S_\hh$
given by 
\begin{equation*}
  \rY := \bigotimes_{i=1}^p \rY^i \,.
\end{equation*}
Then
\begin{equation}
  \label{eq:prod}
  \spec \rY[W(\cH)] = \bigl\{ \l_1 + \cdots + \l_p \,:\, 
  \l_i \in \spec \rY^i[W(\cH_i)] \bigr\} \,.
\end{equation}
\end{proposition}

\Pro\ 
We have
\begin{align*}
  W(\cH) &= \sum_{e\in E(\cH)} e =
  \sum_{i=1}^p \sum_{e\in E(\cH_i)} e = 
  \sum_{i=1}^p W(\cH_i) \\
  &=
  \sum_{i=1}^p \uuno_{S_{\hh_1}} \cdot \cdots \cdot \uuno_{S_{\hh_{i-1}}}
  \cdot W(\cH_i) \cdot  \uuno_{S_{\hh_{i+1}}}
  \cdot \cdots \cdot \uuno_{S_{\hh_p}} \,,
\end{align*}
where $\uuno_G$ stands for the unit element of the group $G$
and of the group algebra $\bC G$.
If $d_i$ is the dimension of the representation $\rY^i$
(which will not be confused, hopefully, with the degree of the vertex $i$),
by the definition of $\rY$ we obtain
\begin{align}
  \label{eq:RW}
  \rY[W(\cH)] &=
  \sum_{i=1}^p \id_{d_1} \otimes \cdots \otimes \id_{d_{i-1}}
  \otimes \rY^i[W(\cH_i)] \otimes  \id_{d_{i+1}}
  \otimes \cdots \otimes \id_{d_p} \,.
\end{align}
Equality \eqref{eq:prod} now follows from a (presumably)
standard argument: 
since $e$ is a transposition, $e^{-1} = e$.
We can assume that representations $\rY^i$ are unitary,
which implies that $\rY^i(e)$ is a Hermitian matrix
and thus $\rY^i[W(\cH_i)]$ is also Hermitian.%
\footnote{Using the Jordan canonical form one can prove
the same result in the general ``non Hermitian'' case.} 
For each $i=1,\ldots,p$, let $(u^{(i)}_j)_{j=1}^{d_i}$
be a basis of $\bC^{d_i}$ consisting of eigenvectors of $\rY^i[W(\cH_i)]$.
The set of all vectors of the form
\begin{equation}
  \label{eq:u}
  u^{(1)}_{j_1} \otimes \cdots \otimes u^{(p)}_{j_p}
\end{equation}
is a basis of $\bC^{\prod_i d_i}$ which consists of
eigenvectors of $\rY[W(\cH)]$. 
Hence the eigenvalues of $\rY[W(\cH)]$
are given \eqref{eq:prod}. 
\qed

\medno
Thanks to identities \eqref{eq:tq}, \eqref{eq:lrp}, Proposition \ref{th:prod}
and Proposition \ref{th:kn} applied to each $K_{\hh_i}$, we have obtained
the following fairly explicit representation
for the eigenvalues of $\rT^\a[W(K_\hh)]$.

\begin{theorem}\label{th:keta}
Let $K_\hh$ be the complete multipartite graph associated
with the partition $\hh=(\hh_1, \ldots, \hh_p)$ of $n$, and let $\rT^\a$ be
one of the (equivalent) irreducible representations of 
$S_n$ corresponding to $\a=(\a_1,\ldots,\a_r)\partit n$.
Then
\begin{align}
  \label{eq:keta}
  \spec \rT^\a[ W(K_\hh) ] =
  \Bigl\{ q_\a - \sum_{i=1}^p q_{\b^i} 
  : \text{$(\b_i)_{i=1}^p$ is $(\a,\hh)$--admissible} 
  \Bigr\}\,.
\end{align}
\end{theorem}

\medno 
We now define, for \textit{arbitrary weak compositions} 
$\b^1, \ldots, \b^p$,
the quantities
\begin{align}
  \label{eq:ba}
  b^\a_{\b^1, \ldots, \b^p} &:= q_\a - \sum_{i=1}^p q_{\b^i} \\
  \label{eq:Ba}
  B^\a_{\hh} &:=
  \max_{(\b^1,\ldots, \b^p)\in \Adm(\a,\hh)} 
  b^\a_{\b^1, \ldots, \b^p} \,.
\end{align}
It follows from Theorem \ref{th:keta} that
\begin{align}
  \label{eq:lmax}
  \lmax\bigl[\a,\, W(K_\hh)\bigr] = B^\a_{\hh} \,.
\end{align}
In order to prove Theorem \ref{th:main} we must show
that
\begin{equation}
  \label{eq:anen}
  \a \ne (n) \quad\impp\quad B^\a_\hh \le B^{(n-1,1)}_\hh \,.
\end{equation}
Our plan at this point is the following: 
we are going to replace the class $\Adm$ with
a different class $\Adm^*$ such that 
\begin{enumerate}[(a)]
\item
the corresponding
maximum
\begin{equation}
  \label{eq:Ba1}
  \bmdf^\a_\hh := 
  \max_{(\b^1,\ldots, \b^p)\in \Adms(\a,\hh)} 
  b^\a_{\b^1, \ldots, \b^p} 
\end{equation}
is easier to evaluate, and, by consequence, we will be able to show
that
\item
\label{it:b}
$\bmdf^\a_\hh$ has a useful monotonicity property with
respect to the \textit{dominance} partial order of partitions. 
A consequence of this monotonicity is that implication \eqref{eq:anen}
holds for the quantities $\bmdf$
\item
\label{it:c}
there is a simple and useful relationship
between the ``true'' quantities $B^\a_\hh$ and their ``fake''
relatives $\bmdf^\a_\hh$, namely:
$B^\a_\hh \le  \bmdf^\a_\hh\;$  and $\;B^\a_\hh =  \bmdf^\a_\hh\;$ if
$\;\a = (n-1,1)$.  
\end{enumerate}
From (\ref{it:b}) and (\ref{it:c}) 
it follows that \eqref{eq:anen} holds,
and hence Theorem \ref{th:main} is proven.

\medno
We start then to describe this alternate class $\Adms$.

\begin{definition}\label{th:adms}
Given two partitions $\a,\hh$ of $n$, with $\hh$ of length $p$,
we denote with $\Adms(\a,\hh)$ the set of all $p$-tuples of
weak compositions $(\g^1, \ldots, \g^p)$ such that
\begin{enumerate}[(i)]
\item $\g^i$ is a weak composition of $\hh_i$
\item $\g^1 + \cdots + \g^p = \a$.
\label{it:adms3}
\end{enumerate}
\end{definition}

\begin{proposition}\label{th:qqq}
Let $\a\partit n$, $\hh=(\hh_1,\ldots,\hh_p)\partit n$. 
\begin{enumerate}[$(1)$]
\item
\label{qqq1}
For each
$(\b^i)_{i=1}^p\in \Adm(\a,\hh)$ 
there exists
$(\g^i)_{i=1}^p\in \Adms(\a,\hh)$ such that
\begin{equation}
  \label{eq:qqq}
  q_{\g^1} + \cdots + q_{\g^p} \le q_{\b^1} + \cdots + q_{\b^p} \,.
\end{equation}
By consequence we have $B^\a_\hh \le  \bmdf^\a_\hh$.
\item
\label{qqq2}
If $\a = (n-1,1)$ and $\hh\ne(1,1,\ldots,1)$, then $B^\a_\hh =  \bmdf^\a_\hh$.
\end{enumerate}
\end{proposition}

\medno
\textit{Proof of (1)}.
The crucial point is the following result that we prove at the
end of this section.

\begin{proposition}\label{th:qq}
Let $n=j+k$ with $j,k$ positive integers and let
$\a\partit n$, $\b\partit j$, $\b'\partit k$ be such that
the Littlewood-Richardson coefficient $c^\a_{\b,\b'}$ is positive.
Then $\b \leC \a$ and
\begin{equation}
  \label{eq:qq}
  q_{\b'} \ge q_{\a-\b} \,.
\end{equation}
\end{proposition}

\medno
Given Proposition \ref{th:qq}, 
we can prove part (\ref{qqq1}) of Proposition \ref{th:qqq} by induction
on $p\ge 2$. If $p=2$, let $\hh=(j,k)$ with $j+k=n$, 
let $(\b,\b')\in \Adm(\a,\hh)$ and define
\begin{equation*}
  (\g,\g') := (\b, \a-\b) \,.
\end{equation*}
The pair $(\g,\g')$ clearly belongs to $\Adms(\a,(j,k))$, and hence
\eqref{eq:qqq} holds thanks to \eqref{eq:qq}.

\smallno
The general case $p\ge 2$ can be proven by induction on $p$.
Assume then that Proposition \ref{th:qqq} holds for $p-1$ and
let $\hh = (\hh_1, \ldots, \hh_p) \partit n$. Define the partitions
\begin{align*}
  \z &= (\z_1, \z_2) := ( \hh_1 + \cdots + \hh_{p-1}, \hh_p ) \vdash n \\
  \hh' &:= (\hh_1, \ldots, \hh_{p-1}) \partit \z_1 \,.
\end{align*}
Since $S_\hh$ is a subgroup of $S_\z$ we have
\begin{equation*}
  [\a] \resR^{S_n}_{S_{\hh}} = 
  \Bigl( [\a] \resR^{S_n}_{S_{\z}} \Bigr) \resR^{S_\z}_{S_{\hh}}   \,.
\end{equation*}
Thus, from the decompositions
\begin{align*}
  [\a] \resR^{S_n}_{S_{\z}} &= 
  \sum_{(\d,\d')\in \Adm(\a,\z)}
    c^\a_{\d,\d'}\; [\d] \otimes [\d'] 
  \\
  [\d] \resR^{S_{\z_1}}_{S_{\hh'}} &= 
  \sum_{(\b^1,\ldots,\b^{p-1})\in \Adm(\d,\hh')}
    c^\d_{\b^1,\ldots,\b^{p-1}}\; [ \b^1] \otimes \cdots \otimes [\b^{p-1}]
\end{align*}
we find
\begin{equation}
  \label{eq:c=sum}
  c^\a_{\b^1,\ldots,\b^p} = \sum_{\d\partit \z_1}
  c^\a_{\d,\b^p} \; c^\d_{\b^1,\ldots, \b^{p-1}} \,.
\end{equation}
If $(\b^1,\ldots,\b^p)\in \Adm(\a,\hh)$, the
quantity in \eqref{eq:c=sum} is positive; 
thus there exists $\d\partit \z_1$
such that both coefficients in the RHS of \eqref{eq:c=sum} are positive.
Pick one such $\d$. 
Since $c^\d_{\b^1,\ldots, \b^{p-1}}>0$, we have,
by induction, that there there exist 
$(\g^1, \ldots, \g^{p-1}) \in \Adms(\d,\hh')$ such that
\begin{equation}
  \label{eq:p-1}
  q_{\g^1} + \cdots + q_{\g^{p-1}} \le q_{\b^1} + \cdots + q_{\b^{p-1}} \,.
\end{equation}
Moreover, $c^\a_{\d,\b^p}>0$, and thus, by Proposition \ref{th:qq}
we have
\begin{equation}
  \label{eq:qbr}
  q_{\b^p} \ge q_{\a - \d} \,.
\end{equation}
Let $\g^p := \a - \d$. From \eqref{eq:p-1}, \eqref{eq:qbr}, and the
fact that
$(\g^1, \ldots, \g^{p-1}) \in \Adms(\d,\hh')$ one can easily conclude
that 
$(\g^1, \ldots, \g^{p}) \in \Adms(\a,\hh)$ and that
inequality \eqref{eq:qqq} holds.

\medno
\textit{Proof of part (\ref{qqq2}) of Proposition \ref{th:qqq}}.
We now show that inequality $B^\a_\hh \le  \bmdf^\a_\hh$
is actually an equality when $\a = (n-1,1)$ and $\hh\ne(1,1,\ldots,1)$.

\smallno
A simple application of the \LR\  rule yields the decomposition
\begin{equation}\begin{split}
  \label{eq:n-1}
  &[n-1,1] \resR^{S_n}_{S_\hh} = 
  (p-1) \; [\hh_1] \otimes \cdots \otimes [\hh_p] \\
  &\quad\oplus
  \bigoplus_{\substack{i=1,\ldots,p \\ \hh_i\ge 2}}
  [\hh_1] \otimes \cdots \otimes [\hh_{i-1}] \otimes[ \hh_i-1,1]
  \otimes [\hh_{i+1}]\otimes\cdots \otimes [\hh_p] \,.
\end{split}\end{equation}
Notice that $[m-1,1]$ has degree $m-1$, while $[\hh_i]$
is the trivial one-dimensional representation, thus
the dimension count is correct in \eqref{eq:n-1}.
From \eqref{eq:n-1} we can read the list of all $(\a,\hh)$--admissible
collections
\begin{equation*}
  \Adm(\a,\hh) := \{ \Psi_0 \} \cup \{ \Psi_i : 1\le i\le p,\ \hh_i\ge 2 \} \,,
\end{equation*}
where
\begin{align}
  \label{eq:qw1}
  \Psi_0 &:= \bigl( (\hh_1), \ldots (\hh_p) \bigr) 
  \\
  \Psi_i &:= \bigl( (\hh_1), \ldots,(\hh_{i-1}), (\hh_i-1,1),(\hh_{i+1}),
  \ldots, (\hh_p) \bigr) \,.
\end{align}
On the other hand, it follows from the definition of $\Adms(\a,\hh)$  that
\begin{equation}
  \label{eq:qw2}
  \Adms(\a,\hh) := \{ \Psi_i : 1\le i\le p  \} \,.
\end{equation}
If all the $\hh_i$'s are greater than 1, then $\Adm^*(\a,\hh)$ is
a subset of $\Adm(\a,\hh)$, and thus the conclusion is trivial.
If some of the $\hh_i$'s are equal to 1, then
in principle we have to worry about the corresponding  
$\Psi_i$, which
belong to $\Adm^*(\a,\hh)$ but not to $\Adm(\a,\hh)$.
But one can easily compute
\begin{align*}
  b^\a_{\Psi_0} &= q_{(n-1,1)} - \sum_{j=1}^p q_{\hh_j} = 
  |E(\cG)| - |V(\cG)|
  \\
  b^\a_{\Psi_i} &= 
  q_{(n-1,1)} - q_{(\hh_i-1,1)} - \sum_{j=1,\, j\ne i}^p q_{\hh_j} =
  |E(\cG)| - |V(\cG)| + \hh_i \,.
\end{align*}
Hence, since, by hypothesis, $\hh_1 = \max_i \hh_i > 1$, we have
\begin{equation*}
  B^{(n-1,1)}_\hh =  \bmdf^{(n-1,1)}_\hh = 
  |E(\cG)| - |V(\cG)| + \hh_1 \,.
  \qquad\qed
\end{equation*}

\begin{remark}\label{th:prov}
Observe that condition $\hh\ne(1,1,\ldots,1)$ is necessary. 
In fact, if $\hh:=(1,1\ldots,1)$, we obtain
\begin{equation*}
  B^{(n-1,1)}_\hh
  = b^\a_{\Psi_0} =
  |E(\cG)| - |V(\cG)| < 
  |E(\cG)| - |V(\cG)| + 1  
   =\bmdf^\a_\hh \,.
\end{equation*}
\end{remark}

\medno
Thanks to condition $\g^1 + \cdots + \g^p = \a$ in the definition
of $\Adm^*(\a,\hh)$, the quantity $b^\a_{\g^1, \ldots, \g^p}$
has a simple expression.

\begin{proposition}\label{th:easy}
If $(\g^1, \ldots, \g^{p}) \in \Adms(\a,\hh)$, we have
\begin{equation*}
  b^\a_{\g^1, \ldots, \g^p} = 
  \sum^p_{\substack{i,j=1 \\  i<j}} \g^i \cdot \g^j \,,
\end{equation*}
where $\g^i \cdot \g^j$ denotes the canonical inner product
$
  \g^i \cdot \g^j := \sum_{k=1}^\oo \g^i_k \, \g^j_k
$.
\end{proposition}

\Pro\ 
The proof is a straightforward computation. Using \eqref{eq:q=} 
and the fact that $\g^1 + \cdots + \g^p = \a$, we find
\begin{align*}
  q_\a - \sum_{i=1}^p q_{\g^i} &=
  \ov2 \sum_{k=1}^\oo \Bigl[ \a_k^2 - \a_k(2k-1) - 
  \sum_{i=1}^p \bigl[ (\g^i_k)^2 - \g^i_k (2k-1) \bigr] \Bigr] \\
  &=
  \ov2 \sum_{k=1}^\oo \Bigl[ \Bigl( \sum_{i=1}^p \g^i_k \Bigr)^2 -
  \sum_{i=1}^p  (\g^i_k)^2  \Bigr] = 
  \sum^p_{\substack{i,j=1 \\  i<j}} \g^i \cdot \g^j \,.
  \qquad\qed
\end{align*}

\medno
In Section \ref{sec:rep} we have defined a ``componentwise'' partial order
$\a\le\b$ in the set of all finite sequences of integers.
We introduce now a weaker
partial order $\leT$, which, following \cite{JaKe},
we call \textit{dominance order}.%
\footnote{One can introduce a third (total) order, namely 
the lexicographic order, but we do not need it in this paper.}
If $\a,\b$ are weak compositions,
we say that $\b$ \textit{dominates} $\a$, and we write $\a\leT\b$,
if%
\begin{equation}
  \label{eq:dom}
  \sum_{i=1}^r (\b_i - \a_i) \ge 0
  \qquad\forall r=1,2,\ldots
\end{equation}
If $\a$ and $\b$ are weak compositions of the \textit{same integer} $n$,
$\b$ dominates $\a$ iff (either $\a=\b$ or)
the Young diagram of $\b$ can be obtained from
the Young diagram of $\a$ by \textit{moving
a certain number of boxes from a lower row to a higher row}.
For instance
\begin{align*}
  \a = \young(\ \ \ \ \ \ ,\ \ \ a,\ \ b,c)
  \leT\ \ 
  \young(\ \ \ \ \ \ ab,\ \ \ ,\ \ c) = \b
\end{align*}
We write $\a \sim \b$ if there exist two integers $j<k$
such that
\begin{equation}
  \label{eq:jk1}
  \text{$\a_i = \b_i$ for all $i\ne j,k$} \qquad\text{and}\qquad
  \b_j - \a_j = \a_k -\b_k = 1\,,
\end{equation}
\textit{i.e.} 
if $\b$ is obtained from $\a$ by removing one box from the right end
of one of its rows and by placing it at the end of a higher row.
Then it is obvious that if $\a\leT\b$,
there exists a finite sequence of ``interpolating'' weak compositions
\begin{equation}
  \label{eq:interp}
  \g_1=\a \leT \g_2 \leT \cdots \leT \g_s = \b
\end{equation}
such that
$\g_i \sim \g_{i+1}$. Less obviously, 
if $\a$ and $\b$ are both partitions of $n$,
then the interpolating sequence $(\g_i)$ can be chosen
in such a way that each $\g_i$ is also a partition of $n$
(see \cite[Theorem 1.4.10]{JaKe}, where a slightly different
notion of $\sim$ is used; his result implies our statement).
The dominance order plays a crucial role in the representation theory
of $S_n$ and, more generally, in combinatorics. We refer the reader
to Section 1.4 of \cite{JaKe}, where several
coimplications of the statement  ``$\a\leT\b$'' are discussed.

\smallno
We are now ready for part (\ref{it:b}) of our ``plan''
outlined above.

\begin{proposition}\label{th:mono}
Let 
$\a,\b$ be two partitions of $n$ with $\a\leT\b$.
\begin{enumerate}[$(1)$]
\item
\label{mono1}
$q_\a \le q_\b$.
\item
\label{mono2}
$\bmdf^\a_\hh \le \bmdf^\b_\hh$ for any  $\hh\partit n$.
\end{enumerate}
\end{proposition}

\noindent
\textit{End of proof of Theorem \ref{th:main}}.
Before proving Proposition \ref{th:mono} we complete the proof
of Theorem \ref{th:main}. 
We start with the simple observation that
\begin{equation}\begin{split}
  \label{eq:simp}
  &\text{\it If $\a$ is a nontrivial partition of $n$,}\\[-1mm]
  &\text{\it \iie\  $\a \ne (n)$,  then $\a \leT (n-1,1)$.}
\end{split}\end{equation}

\smallno
If $\hh=(1,1,\ldots,1)$, $\cG$ is the
complete graph $K_n$. 
This case is well known \cite{DiSh}.
Anyway, the proof goes as follows:
Proposition \ref{th:kn} says that $\rT^\a[W(K_n)]$
has a unique eigenvalue $q_\a$ of multiplicity equal
to the dimension of the representation $\a$; in particular,
$\lmax[\a,\, W(K_n)] = q_\a$. 
Statement (\ref{mono1})
of Proposition \ref{th:mono} yields $\lmax[\a,\, W(K_n)] \le
\lmax[(n-1,1),\, W(K_n)]$.

\smallno
Case $\hh\ne(1,1,\ldots,1)$.
Thanks to \eqref{eq:lmax}, \eqref{eq:simp}, Proposition \ref{th:qqq} and Proposition 
\ref{th:mono} we have
\begin{align*}
  &\lmax\bigl[\a,\, W(K_\hh)\bigr] = 
  B^\a_\hh \le \bmdf^\a_\hh \le \bmdf^{(n-1,1)}_\hh \\
  &\quad= 
  B^{(n-1,1)}_\hh =
  \lmax\bigl[(n-1,1),\, W(K_\hh)\bigr] 
  \,.
  \qquad\qed
\end{align*}

\medno
\textit{Proof of Proposition \ref{th:mono}}.
Thanks to 
the existence of the interpolating sequence \eqref{eq:interp},
we can assume that $\a \sim\b$.
Thus there exist two integers $j<k$ such that
\begin{equation}
  \label{eq:jk3}
  \text{$\a_i = \b_i$ for all $i\ne j,k$} \qquad\text{and}\qquad
  \b_j - \a_j = \a_k -\b_k = 1\,.
\end{equation}
Statement (\ref{mono1}) was proven in
\cite[Lemma 10]{DiSh} and it is a simple computation.
In fact from \eqref{eq:q=}, \eqref{eq:jk3}, and the fact that 
$\a$ is a partition it follows that
\begin{equation*}
  q_\b - q_\a = (\a_j - \a_k) + (k - j + 1) \ge k-j+1 \ge 2 \,.
\end{equation*}
For the proof of (\ref{mono2}),
let $\hh=(\hh_1,\ldots,\hh_p)\partit n$. 
We now show that for each $(\g^i)_{i=1}^p \in \Adms(\a,\hh)$ 
there exists $(\d^i)_{i=1}^p\in \Adms(\b,\hh)$
such that 
\begin{equation}
  \label{eq:babb}
  b^\a_{\g^1, \ldots, \g^p} \le 
  b^\b_{\d^1, \ldots, \d^p} \,.
\end{equation}
Let then $(\g^1, \ldots, \g^p) \in \Adms(\a,\hh)$ and,
for a positive integer $\ell$, which will be determined
at the end,
define the $p$-tuple $(\d^i)_{i=1}^p$ as follows: 
\begin{enumerate}[(a)]
\item if $i\ne \ell$ we simply let $\d^i = \g^i$, while
\item $\d^\ell$ is obtained from $\g^\ell$ by moving
one box from row $k$ to row $j$, \iie\ 
\begin{align}
  \label{eq:jk2}
  &\text{$\d^\ell_i = \g^\ell_i$ for all $i\ne j,k$} 
  &\d^\ell_j &= \g^\ell_j+1  &  \d^\ell_k &= \g^\ell_k -1 \,.
\end{align}
\end{enumerate}
It is clear that $(\d^i)_{i=1}^p \in \Adms(\b,\hh)$ \textit{unless}
$\g^\ell_k=0$, but we will worry about this later.
Using Proposition \ref{th:easy} we obtain
\begin{equation}\begin{split}
  \label{eq:bb}
  &\D_\ell b := b^\b_{\d^1, \ldots, \d^p} - 
  b^\a_{\g^1, \ldots, \g^p}  =
  \sum^p_{i,m=1 ,\,  i<m} 
   \[ \d^i \cdot \d^m -
      \g^i \cdot \g^m \] \\
  &\quad
  = \bigl( \d^\ell - \g^\ell\bigr) \cdot \sum^p_{i=1,\, i\ne \ell} \g^i =
  \sum^p_{i=1,\, i\ne \ell} \bigl( \g^i_j - \g^i_k \bigr) 
  \\
  &\quad=
  (\a_j-\a_k) - (\g^\ell_j - \g^\ell_k) \,.
\end{split}\end{equation}
From \eqref{eq:jk3} it follows that $\a_k \ge 1$, so the set
of ``legal'' values of $\ell$, that is
those such that $\g^\ell_k >0$, is nonempty.
Summing over all these values of $\ell$ and keeping in mind
that $\a_j \ge\a_k$, we have
\begin{align*}
  &\sum_{\ell:\, \g^\ell_k > 0} \D_\ell b =
  |\{ \ell : \g^\ell_k>0 \}| \: (\a_j - \a_k) - 
  \sum_{\ell:\, \g^\ell_k > 0} \g^\ell_j + \a_k
  \\
  &\quad\ge
  \a_j - \sum_{\ell:\, \g^\ell_k > 0} \g^\ell_j 
  = \sum_{\ell:\, \g^\ell_k = 0} \g^\ell_j 
  \ge 0\,.
\end{align*}
Hence there exists at least one value of $\ell$ such that
$\g^\ell_k>0$ and $\D_\ell b \ge 0$.
If we then define $(\d_i)_{i=1}^p$ using this value of $\ell$, 
we get \eqref{eq:babb}.
\qed

\bigno
\textit{Proof of Proposition \ref{th:qq}}

\noindent
Let $n=j+k$, $\a\partit n$, $\b\partit j$, $\g\partit k$,
and assume that
the Littlewood-Richardson coefficient $c^\a_{\b,\g}$ is positive.
The fact that $\b \leC \a$, \iie\  that $\b_i \le \a_i$ for each $i$,
trivially follows from the \LR\ rule which we discuss in Section
\ref{sec:lr}.
We must prove that $q_{\g} \ge q_{\a-\b}$.
The difference $\a-\b$ is a weak composition of $k$ but
not necessarily a partition.
We can nevertheless obtain a partition by sorting the entries of $\a-\b$
in nonincreasing order (and dropping the trailing zeros). 
We denote this partition with
$\srt(\a-\b)$. So, if
\begin{equation*}
  \a - \b = \d = (\d_1, \ldots, \d_r)
\end{equation*}
we have
\begin{align}
  \srt(\d) &= (\d_{\pi(1)}, \ldots, \d_{\pi(r)})
  &
  \text{with}
  &
  &\d_{\pi(1)}\ge \cdots \ge \d_{\pi(r)} \,,
\end{align}
where $\pi$ is a suitable permutation in $S_r$.
We start by observing that this sorting procedure does not decrease
the quantity $q$.

\begin{proposition}\label{th:qinc}
Let $\a,\b$ be two weak compositions of $n$. 
Assume that $\b$ is obtained from $\a$ by switching two elements 
$\a_j, \a_k$ with $j<k$, \iie\ 
\begin{equation*}
  \b = (\a_1, \ldots, \a_{j-1},\a_k,\a_{j+1},\ldots, \a_{k-1},\a_j,
  \a_{k+1},\ldots) \,.
\end{equation*}
Then%
\footnote{Remember that we have \textit{defined}
$q_\a$ for weak compositions by formula \eqref{eq:q=}.}
$\a_j < \a_k$ if and only if  $q_\a < q_\b$.
By consequence $q_{\srt(\a)} \ge q_\a$ for any weak composition $\a$.
\end{proposition}

\Pro\
From \eqref{eq:q=} we get
\begin{align*}
  q_\b - q_\a &= 
  \ov2 \bigl[ \a_k(\a_k - 2j + 1)   - \a_j (\a_j - 2j+1)  
  \\&
  \quad +
          \a_j(\a_j - 2k+1) - \a_k (\a_k - 2k+1) \bigr]
  \\&
  = (\a_j-\a_k) \, (j-k)  \,.
\end{align*}
The proposition follows, since $j<k$.
\qed

\bigno
The central ingredient of the proof of Proposition \ref{th:qq}
is the following property
of the Richardson-Littlewood coefficients, which we prove
in Section \ref{sec:lr}, and which identifies 
the partition $\g$ with ``minimal content''.

\begin{lemma}\label{th:lrmin}
Let $\a\partit n$ and  $\b\partit j$, with $j<n$ and $\b \leC \a$.
Let $\hg := \srt(\a-\b)$, so that $\hg$ is a partition of $k=n-j$.
\begin{enumerate}[$(1)$]
\item 
The \LR\ coefficient $c^\a_{\b,\hg}$ is positive.
\item
If $\g$ is a partition of $k$ such that $c^\a_{\b,\g}>0$,
then $\g \geT \hg$.
\end{enumerate}
\end{lemma}

\medno
Given Lemma \ref{th:lrmin}, the proof of Proposition \ref{th:qq}
readily follows. In fact, from $\g\geT\hg$, statement (\ref{mono1})
of Proposition \ref{th:mono} and
Proposition \ref{th:qinc} we obtain
\begin{equation*}
  q_\g \ge q_{\hg} = q_{\srt(\a-\b)} \ge q_{\a-\b} \,. 
  \qquad\qed
\end{equation*}

\begin{remark}
Theorem \ref{th:main} states that, if $\cG$ is complete
multipartite, then
\begin{equation*}
  \a \leT (n-1,1) \quad\impp\quad
  \lmax\bigl[ \a,\, W(\cG) \bigr]  \le
  \lmax\bigl[ (n-1,1),\, W(\cG) \bigr] \,.
\end{equation*}
One may wonder wheth\-er the following
``strict version'' of this implication holds:
\begin{equation}
  \label{eq:af2}
  \a \triangleleft (n-1,1) \quad\impp\quad
  \lmax\bigl[ \a,\, W(\cG) \bigr]  <
  \lmax\bigl[ (n-1,1),\, W(\cG) \bigr] \,.
\end{equation}
This would imply that the first nontrivial eigenvalues
of $\D_\cG$ and $\D_{\Cay(\cG)}$ are not only equal but
also have 
\textit{the same multiplicity}.
But \eqref{eq:af2} is false. Take, in fact, $K_{2,2}$.
It is easy to verify that
\begin{equation*}
  \lmax\bigl( (2,2),\, W[K_{2,2}] \bigr) =
  \lmax\bigl( (3,1),\, W[K_{2,2}] \bigr) = 2 \,.
\end{equation*}
\end{remark}

\secskip
\section{\LR\ tableaux with minimal content}
\label{sec:lr}

\subsection{The \LR\ rule}

Since we are going to deal with the mechanisms of the \LR\ rule,
it might be a good idea to briefly describe how it works.
We need to define \textit{skew tableaux}, 
\textit{semistandard tableaux}, \textit{lattice permutations}
and \textit{content} of a tableau. Impatient readers who
are not acquainted with this rule, and who prefer
worked examples to abstract definitions,
might try staring at Table \ref{tab:lr}
for a couple of minutes.%
\footnote{It didn't work for the author though!}

Let as usual $j,k$ be two positive integers with $n=j+k$.
Given $\a \partit n$ and $\b\partit j$ with $\b \leC \a$, 
the \textit{skew diagram}
of shape $\a/\b$ is the set of boxes obtained by erasing
in the Young diagram $\a$ all boxes which also appear in $\b$.
Let, for example, $\a=(7,6,3,1)$ and $\b=(5,2,1)$. Let us
draw $\a$ and cross all boxes which belong to $\b$. The
skew diagram of shape $\a/\b$ is the set of all uncrossed boxes
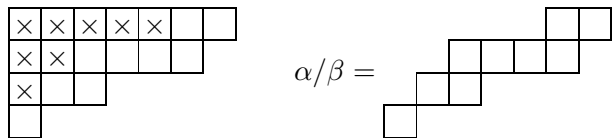
\begin{equation*}
  \let\x\times
  \young(\x\x\x\x\x\  \ ,\x\x\ \ \ \ ,\x\ \ ,\ ) \qquad
  \a/\b = \young(:::::\ \ ,::\ \ \ \ ,:\ \ ,\ )
\end{equation*}
A \textit{skew tableau} is a skew diagram with a positive
integer placed in each box. A skew tableau is called
\textit{semistandard} if each row is nondecreasing and
each column is strictly increasing.
A \textit{lattice permutation} is a finite sequence of
positive integers
$\o = (\o_1, \o_2, \ldots, \o_s)$ such that for each $k=1,\ldots,s$
the number of times that any given integer $i$ appears
in the initial subsequence $(\o_1,\ldots\o_k)$ cannot 
exceed the number of times that $i-1$ appears in the same subsequence.
In other words
\begin{equation*}
  | \{ j\le k : \o_j = a \} | \ge | \{ j\le k : \o_j = b \} |
  \qquad \forall k=1,\dots,s
  \qquad \forall a < b \,.
\end{equation*}
For instance, if we let
\begin{equation*}
  \o := (1,1,2,1,3,2,2,3,1,1,2)
  \qquad
  \o' := (1,1,2,1,3,2,3,\uu 3,1,1,2)
\end{equation*}
then $\o$ is a lattice permutation, but $\o'$ is not, since,
when we arrive at the underlined digit $\uu 3$, we realize 
that we have encountered along the way more $3$'s than  $2$'s.


\begin{table}
\begin{center}
\let\x=\times
\renewcommand\arraystretch{1.0}
\caption{LR tableaux of shape $\a/\b$ and content $\g$}
\label{tab:lr}
\begin{tabular}{|c|c|c|c|} \hline
\multicolumn{4}{|c|}{$\vphantom{\Bigl[}%
\a=(6,5,3,1) \qquad
\b=(5,2,1)$} \\ \hline
  $\vphantom{\Bigl[}$
  LR tableau & $\g$
&
  LR tableau & $\g$
\\ \hline
\ &&&\\
       \yuung{:::::1,::111,:12,1}
     & (6,1)
     & \yuung{:::::1,::111,:12,2}
     & (5,2)
\\
\ &&&\\
       \yuung{:::::1,::111,:22,1}
     & (5,2)
     & \yuung{:::::1,::112,:12,1}
     & (5,2)
\\
\ &&&\\
       \yuung{:::::1,::111,:12,3}
     & (5,1,1)
     & \yuung{:::::1,::112,:13,1}
     & (5,1,1)
\\
\ &&&\\
       \yuung{:::::1,::111,:22,2}
     & (4,3)
     & \yuung{:::::1,::112,:12,2}
     & (4,3)
\\
\ &&&\\
       \yuung{:::::1,::112,:22,1}
     & (4,3)
     & \yuung{:::::1,::111,:22,3}
     & (4,2,1)
\\
\ &&&\\
       \yuung{:::::1,::112,:12,3}
     & (4,2,1)
     & \yuung{:::::1,::112,:13,2}
     & (4,2,1)
\\
\ &&&\\
       \yuung{:::::1,::112,:23,1}
     & (4,2,1)
     & \yuung{:::::1,::112,:13,4}
     & (4,1,1,1)
\\
\ &&&\\
       \yuung{:::::1,::112,:22,3}
     & (3,3,1)
     & \yuung{:::::1,::112,:23,2}
     & (3,3,1)
\\
\ &&&\\
       \yuung{:::::1,::112,:23,3}
     & (3,2,2)
     & \yuung{:::::1,::112,:23,4}
     & (3,2,1,1)
\\\ &&&\\
\hline
\end{tabular}
\end{center}
\end{table}

The \textit{content} of a finite sequence of positive integers 
$\o =(\o_1, \ldots, \o_s)$ is the sequence
$\g=(\g_1,\g_2,\ldots)$, where
$\g_i$ is the number of times the integer $i$
appears in $\o$. 
The \textit{content} of a (skew) Young tableau $t$
is the content of the sequence of all the integers which appear
in $t$, listed in any (obviously arbitrary) order.
So, for instance,
\begin{equation*}
  \young(::::1356,:::2131,:12)
  \quad\text{has content}\quad  (4,2,2,0,1,1)
\end{equation*}
We are finally able to state the \LR\ rule.
For a proof see, for instance, \cite{JaKe} or \cite{Sag}.

\begin{LRR}
The coefficient $c^\a_{\b,\g}$ which appears in \eqref{eq:lr} 
is equal to the number of semistandard
skew tableaux of shape $\a/\b$ and content $\g$, which yield
lattice permutations when we read their entries from \uu{right to left}
and downward.%
We will call these tableaux \textit{LR tableaux}.
\end{LRR}

\bigbreak

\noindent
Table \ref{tab:lr} shows all LR tableaux of shape $\a/\b$ with content $\g$,
for $\a=(6,5,3,1)$ and $\b=(5,2,1)$. Notice that the minimal $\g$
is the last entry of the table, $\g=(3,2,1,1) = \srt(\a-\b)$,
in accord with Lemma \ref{th:lrmin}.

\subsection{Proof of Lemma \ref{th:lrmin}}

Let $t$ be an LR tableau of shape $\a/\b$ and let $\d=\a-\b$.
Since $\b\leC\a$, in general $\d$ is a \textit{weak} 
composition of $k:=|\a|-|\b|$,
so it can have (nontrailing) zeros. 
In the next remark we get rid of these zeros, reducing to the
case $\d \scomp k$. This is really irrelevant, as we will see,
but it simplifies some statements.
For this purpose we introduce the notation 
\begin{equation*}
  \b\leCC\a \quad\coimpdef\quad
  \text{$\b_i < \a_i$ for all $i$ such that $\a_i>0$.}
\end{equation*}
so that $\b\leCC\a$ implies that $\a-\b$ is a composition.

\begin{remark}\label{th:redu}
(\textit{Reduction to the case $\b\leCC\a$}).
Assume that the Young diagrams of $\a$ and $\b$ have some
corresponding rows of equal length $\a_i = \b_i$. 
Let $\ol \a$ and $\ol \b$
be the Young diagrams obtained by eliminating 
in $\a$ and $\b$ all corresponding rows of equal length.
Then it is obvious that the number of LR tableaux of
shape $\a/\b$ and content $\g$ is not affected by the
simultaneous replacements of $\a$ with $\ol\a$
and $\b$ with $\ol\b$. Since, on the other hand,
$\srt(\a-\b) = \srt(\ol\a-\ol\b)$, it follows
that if Lemma \ref{th:lrmin} holds when $\b\leCC\a$,
it also holds for $\b\leC\a$. 
\end{remark}

\medno
It is our intention
to concatenate the rows of an LR tableau $t$ into a single sequence
of positive integers, of length $|\d|$,
that we denote by $\th(t)$. Rather than keep struggling
with our primordial instinct
to read (\textit{and think}) from left to right,%
\footnote{We apologize to native right-to-left thinkers.}
we first flip each
row of $t$, and then concatenate the rows from top to
bottom. So, for instance, if
\begin{equation}
  \label{eq:abcd}
  \a=(7,6,4,3) \qquad \b=(4,2,1,1)\qquad 
  \g=(6,3,2,1) \qquad \d=(3,4,3,2) \,,
\end{equation}
an example of an LR tableaux $t$ of shape $\a/\b$ and content $\g$,
and relative sequence $\th(t)$ is
\begin{align}
  \label{eq:t}
  t &= \young(::::111,::1122,:123,:34) \\[3mm]
  \th(t) &= (\, 
  \ubr{1,1,1}{N^\d_1} \,,\,  
  \ubr{2,2,1,1}{N^\d_2} \,,\, 
  \ubr{3,2,1}{N^\d_3} \,,\, 
  \ubr{4,3}{N^\d_4} \,)
\end{align}
We observe that $\th(t)$ is nonincreasing in each 
subinterval corresponding a row of $t$. 


\begin{definition}\label{th:dd}
Given $\d=(\d_1,\ldots, \d_r)\scomp n$, 
we say that a sequence of $n$ positive integers
is \textit{$\d$--nonincreasing} if it is nonincreasing
in each subinterval $N^\d_i$ (recall definition \eqref{eq:N}).
For future purposes, it is convenient
to stretch a little bit (in an obvious way) this definition, in order 
to accommodate also sequences shorter than $|\d|$.
If $m<n = |\d|$, a sequence of $m$ positive integers is called
\textit{$\d$--nonincreasing} if the sequence
of length $n$ obtained by adding $n-m$ $1$'s at the
end is $\d$--nonincreasing.
\end{definition}

\smallno
We denote with $\O_\d$ the set of all sequences 
$\o=(\o_1, \ldots, \o_{|\d|})$ of positive integers such that
\begin{enumerate}[(i)]
\item $\o$ is a lattice permutation
\item $\o$ is $\d$--nonincreasing.
\end{enumerate}

\smallno
We have then:

\begin{proposition}\label{th:O}
Given $\a\partit n$ and $\b \partit j$ with $\b\leCC \a$,
the mapping $t\mapsto \th(t)$ is an injection
from the set of all LR tableaux of shape $\a/\b$
into $\O_{\a-\b}$.
\end{proposition}

\Pro\ 
It is obvious from the construction of $\th(t)$
and from the definition of LR tableau.
\qed

\smallno
\begin{remark}\label{th:bad}
Notice that in general, given $\a$ and $\b$, 
and letting $\d:=\a-\b$,
the mapping
\begin{equation*}
  \th : \{ \text{LR tableaux of shape $\a/\b$} \} \mapsto \O_\d 
\end{equation*}
is not surjective. There will be in fact elements
$\o \in \O_{\d}$ that, when ``disassembled back and flipped''
in order to restore a tableau,
are going to produce a tableau $t$ which fails to satisfy the
requirement of having increasing columns.
In particular $\o = (1,1,\ldots,1)$ is always an element
of $\O_{\d}$, but it is not in the image of $\th$
unless $\b_i \ge \a_{i+1}$ for all $i$, which
means that the rows of skew tableaux of shape $\a/\b$
are shifted in such a way that they do not ``overlap'',
like, for instance
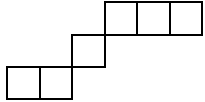
\begin{equation*}
  \young(:::\ \ \ ,::\ ,\ \ ) 
\end{equation*}
So, given $\a$ and $\b$, the set $\O_\d$ contains ``good''
sequences as well as ``extraneous'' sequences
which do not correspond to any LR tableau of shape $\a/\b$.
On the other hand the \textit{advantage} of dealing
with $\O_\d$ is that this set depends
on the pair $(\a,\b)$ only through
their difference $\a-\b$. In the end, this works because,
quite remarkably, the solution to the problem
of finding the minimal
$\g$ such that $c^\a_{\b,\g}>0$
does indeed depend only on $\a-\b$!
\end{remark}

\medno
We now describe a simple algorithm for
constructing, for a given $\d\scomp k$, a particular 
sequence $\ho \in \O_\d$ which we has the property
of having \textit{minimal content} with respect
to all other $\o\in \O_\d$. 
Then we show that the content of $\ho$ is $\srt(\d)$.
Finally we prove that
$\ho$ is a \textit{good} sequence, meaning
that there exists an LR tableau $t$ such that $\ho = \th(t)$,
completing in this way the proof of Lemma \ref{th:lrmin}.

\smallno
Let then $\d$ be a composition of $k$. The 
sequence $\ho$ in $\O_\d$ is defined recursively as follows:
\begin{enumerate}[(a)]
\item start with $\ho_1 = 1$ (mandatory)
\item given $(\ho_1, \ldots, \ho_{i-1})$, choose
$\ho_i$ as the largest integer such that the sequence
$(\ho_1, \ldots, \ho_i)$ satisfies 
(i) and (ii) in the definition
of $\O_\d$.
\end{enumerate}
We call $\ho$ the \textit{minimal sequence} in $\O_\d$.%
\footnote{The reader may rightfully feel that we
are testing his patience, by calling \textit{minimal}
a sequence where we pick at each step the largest possible
value. But, again, $\ho$ has minimal \textit{content}.
In case someone were wondering, $\ho$ is \textit{not} maximal
in $\O_\d$, in general,  with respect to the componentwise partial order.%
}
We observe that this algorithm can never get stuck. 
Let in fact $X_i$ be the 
set of possible choices for $\ho_i$, \iie\  the set of integers
$s$ such that $(\ho_1, \ldots, \ho_{i-1}, s)$
satisfies (i) and (ii) in the definition
of $\O_\d$. Then $X_i$
is nonempty
($1$ is always a possible choice) and bounded (a trivial
bound for being a lattice permutation is $s=\ho_i\le i$), thus
the maximum of $X_i$ exists.

\smallno
For example, if $\d = (3,4,2,4)$, the minimal sequence is given by
\begin{equation}
  \label{eq:ho}
  \ho = 
  (\, 
  \ubr{1,1,1}{N^\d_1} \,,\,  
  \ubr{2,2,2,1}{N^\d_2} \,,\, 
  \ubr{3,3}{N^\d_3} \,,\, 
  \ubr{4,4,3,2}{N^\d_4} 
  \,) \,.
\end{equation}
Notice that, at least in this case, the content of $\ho$
is $(4,4,3,2) = \srt(\d)$.

\begin{proposition}\label{th:ho}
Let $\ho$ be the minimal sequence in $\O_\d$.
\begin{enumerate}[$(1)$]
\item
If $z_\ell$ is the first integer of the subinterval $N^\d_\ell$,
we have $\ho_{z_\ell}=\ell$.
\item
$\max \{ \ho_i : i \in N^\d_{\ell} \} = \ell$.
\item
If $p:=\ho_i > \ho_{i+1} =:m$, then
$| \{ j\le i : \ho_j = m \} | = | \{ j\le i : \ho_j = p \} |$.
\end{enumerate}
\end{proposition}

\Pro\ 
Let $z_\ell = \d_1 + \cdots + \d_{\ell-1} + 1$ be
the first integer in $N^\d_\ell$.
For $\ell=1$ we have $z_\ell=1$ and, necessarily, $\ho_1=1$.
Assume then  that $\ho_{z_j}=j$ for $j=1, \ldots,\ell-1$.
Since $\ho$ is $\d$--nonincreasing, we have
\begin{equation}
  \label{eq:maxo}
  \max\{ \ho_i : i \in N^\d_1 \cup \cdots \cup N^\d_{\ell-1} \} = 
  \ho_{z_{\ell-1}} = \ell-1 \,.
\end{equation}
Imagine to have constructed $(\ho_1, \ldots, \ho_{z_\ell-1})$.
Since $z_\ell$ is at the beginning
of $N^\d_{\ell}$, for picking $\ho_{z_\ell}$, 
we do not worry about $\d$--nonincreasing.
On the other hand, from \eqref{eq:maxo}, it follows 
that each integer from $1$ to $\ell-1$ has already
appeared in $\ho$ at least once, while no integer greater than or equal
to $\ell$ has yet appeared. 
By consequence, $\ho_{z_\ell} = \ell$.
This proves part (1) by induction.
Statement (2) follows from (1) and from the fact that $\ho$ is 
$\d$--nonincreasing.

\smallno
Assume $p = \ho_{i} > \ho_{i+1} = m$. This means,
by the definition of $\ho$, that none of the integers
from $m+1$ to $p$ can be picked as $\ho_{i+1}$, since
that would violate the lattice permutation property.
Hence we must have
\begin{align*}
  &| \{ j\le i : \ho_j = m-1 \} | >
  | \{ j\le i : \ho_j = m \} | \\
  &\quad
  = | \{ j\le i : \ho_j = m+1 \} | =
  \cdots = | \{ j\le i : \ho_j = p \} | \,.
  \qquad\qed
\end{align*}

\smallno

\begin{proposition}\label{th:mincont}
Let 
$\o \in \O_\d$ for some $\d \scomp n$, 
and let $\g$ be the content of $\o$. 
Let also $\hg$ be the content of the minimal sequence  $\ho$.
Then
\begin{equation*}
  \g \geT \hg = \srt(\d) \,.
\end{equation*}
\end{proposition}

\Pro\ 
We first observe that, given two arbitrary finite sequences $\o, \ol\o$
of positive integers with respective contents $\b$ and $\ol\b$,
we have
\begin{equation}
  \label{eq:oo}
  \o\le\ol\o 
  \quad \coimpp \quad
  \b \geT \ol\b \,.
\end{equation}
This simply follows from the definition \eqref{eq:dom} of dominance
and from the identity $\sum_{i=1}^s \b_i = |\{ j: \o_j \le s \}|$.

It is well known 
\cite[1.4.11]{JaKe} that the dominance order is reversed on
conjugate partitions \eqref{eq:dual}, that is
\begin{equation}
  \label{eq:abba}
  \a \leT \b \quad\coimpp\quad \a' \geT \b' \,,
\end{equation}
and thus we are going to prove that $\g' \leT \hg'$.
Given $\o=(\o_1, \ldots, \o_n)\in \O_\d$, we
define the \textit{running multiplicity} of $\o$ as the 
sequence $\nu=(\n_1, \ldots, \n_n)$ where
\begin{equation*}
  \nu_i := | \{ j\le i : \o_j = \o_i \}| \,.
\end{equation*}
In other words $\nu_i$ is the number of times that $\o_i$ has appeared in
the sequence ``up to that point''.
Proposition \ref{th:mincont} is a consequence of the
following Lemma.

\begin{lemma}\label{th:abcd}
\begin{enumerate}[(A)]
\item 
The content of $\nu$ is $\g'$, that is
$\g'_i = |\{ j : \nu_j = i \}|$. In other words,
\textit{the content of the running multiplicity
of $\o$ is the dual of the content of $\o$.} 
\item
The sequence $\nu$ is strictly increasing on each
subinterval of $\{1,\ldots,n\}$ on which $\o$ is
nonincreasing. 
In particular, since $\o$ is $\d$--nonincreasing, 
$\nu$ is strictly increasing
on each subinterval $N^\d_i$.
\item
If $\hn$ is the running multiplicity of $\ho$ we have, 
denoting the length of $\d$ with $r$,
\begin{equation}
  \label{eq:hn}
  \hn = 
  (\, 
  \ubr{1,2,3,\dots,\d_1}{N^\d_1} \,,\,  
  \ubr{1,2,3,\dots,\d_2}{N^\d_2} \,,\,  
  \ldots \,,\,
  \ubr{1,2,3,\dots,\d_r}{N^\d_r}
  \,)
\end{equation}
\item
$\hg=\srt(\d)$.
\end{enumerate}
\end{lemma}

\smallno
From (A), (B), (C), and (D) Proposition \ref{th:mincont}
follows. In fact,
\begin{align*}
  (B) +(C)  \quad&\impp\quad \hn \le \nu  
  \\
  \quad&\impp\quad \hg' \geT \g'
  &\text{[thanks to (A) and \eqref{eq:oo}]}
  \\
  \quad&\impp\quad \g \geT \hg = \srt(\d)
  &\text{[thanks to \eqref{eq:abba} and (D)].}
\end{align*}
We are then left with the proof Lemma \ref{th:abcd}.

\smallno
Proof of (A). 
$\g_i\ge m$ means that $i$ appears at least $m$ times in $\o$.
Hence there exists a positive integer $k_i\le n$ which marks
the $m^{\text{th}}$ appearance of $i$ in $\o$. Clearly $\nu_{k_i} = m$.
Viceversa if $\nu_k = m$ the integer $\o_k$ appears $m$
times in $\{\o_1, \ldots, \o_k\}$, thus
we have $\g_{\o_k} \ge m$.
We have then, for fixed $\o$ and $m$, a bijection
\begin{equation*}
  \{ i : \g_i \ge m \} \ni i \otto k_i \in \{ k : \nu_k = m \} \,,
\end{equation*}
which implies
\begin{equation*}
  \g'_m = |\{ i : \g_i \ge m \}| = |\{ k : \nu_k = m \}| \,,
\end{equation*}
\iie\  that $\g'$ is the content of $\nu$.

\smallno
Proof of (B). 
If $\o_{i+1} \le \o_{i}$, since $\o$ is a lattice permutation,
we have
\begin{align*}
  \n_{i+1} &=
  |\{ j \le i+1 : \o_j = \o_{i+1} \}|  =
  |\{ j \le i : \o_j = \o_{i+1} \}| + 1\\
  &\ge
  |\{ j \le i : \o_j = \o_i \}| + 1
  = \n_i +1 \,.
\end{align*}

\smallno
Proof of (C). 
Let $z_\ell$ be the first integer in $N^\d_\ell$.
Proposition \ref{th:ho} 
implies that the integer $\ell$ appears for the first
time in $\ho$ at position $z_\ell$, and hence $\hn_{z_\ell}=1$.

\smallno
In order to prove \eqref{eq:hn} we show that, 
within each subinterval $N^\d_\ell$,
we have $\hn_{i+1} = \hn_{i}+1$. 
If $\ho_{i+1} = \ho_i$,
this is a trivial 
consequence of the definition of $\hn$.
On the other hand, if $p:=\ho_{i} > \ho_{i+1}=:m$, thanks
to Proposition \ref{th:ho}, we have
\begin{align*}
  \hn_{i+1} &= | \{ j\le i+1 : \o_j = m \} | 
  = | \{ j\le i : \o_j = m \} | + 1 \\
  &= | \{ j\le i : \o_j = p \} | + 1 = \hn_i + 1 \,.
\end{align*}

\smallno
Proof of (D).
Let $\eps := \srt(\d)$, so that $\eps$ is a partition of $|\d|$.
From statement (A) applied to $\hn$,
from the explicit expression of $\hn$ given by \eqref{eq:hn},
and from the definition of duality \eqref{eq:dual}, 
it follows that the content of $\hn$ is
given by 
\begin{equation*}
  \hg'_i 
  = |\{ j : \hn_j = i \} |
  = |\{ j : \d_j \ge i \} |= |\{ j : \eps_j \ge i \} |= \eps'_i \,.
\end{equation*}
But $\g$ and $\eps$ are both partitions, so we can
take the dual and get (D). Lemma \ref{th:abcd}
and Proposition \ref{th:mincont} are thus proven.
\qed

\bigno
In order to complete the proof of Lemma \ref{th:lrmin},
one last step is required, namely to check that
the minimal sequence $\ho$ in $\O_\d$
corresponds to an actual LR tableaux.

\begin{proposition}\label{th:good}
Let $\d$ be a composition of $k$, 
and let $\ho$ be the minimal sequence in $\O_\d$
with content $\hg = \srt(\d)$.
If $\a$ is a partition of $n$ and $\b$ is a partition of $n-k$
such that
$\a - \b = \d$, then there exists an LR tableau $t$ of shape $\a/\b$
and content $\hg$ such that $\th(t) = \ho$.
\end{proposition}

\Pro\ 
Let $z_\ell$ be the first integer
in the interval $N^\d_\ell$. Hence
\begin{equation}
  \label{eq:ho1}
  \ho = 
  (\, 
  \ubr{\ho_{z_1} \dots, \ho_{\d_1}}{N^\d_1} \,,\,  
  \ubr{\ho_{z_2},\dots, \ho_{z_2+\d_2-1}}{N^\d_2} \,,\,  
  \ldots \,,\,
  \ubr{\ho_{z_r},\dots, \ho_{z_r+\d_r-1}}{N^\d_r} 
  \,)
\end{equation}
In order to reconstruct a tableau $t$ from $\ho$ we proceed in three steps.
\begin{enumerate}[(s1)]
\item
We arrange the restrictions of $\ho$ to each subinterval $N^\d_\ell$
as the rows of an (improper)
tableau with a \textit{left aligned} border. 
We denote this tableau with $\ds\Tt$.
\item
For each $\ell=2, \ldots,r$, 
we shift the $\ell^{\rm th}$ row of $\Tt$ 
by $\a_1 - \a_\ell$ positions to the right, and denote this tableau 
with $\tshi$.
\item
Finally we flip $\tshi$ horizontally and get the tableau $t$
which is a skew tableau of shape $\a/\b$.
\end{enumerate}
Consider for example the sequence $\ho$ in \eqref{eq:ho}.
If $\a = (8,7,5,4)$,  transformations (s1), (s2) and (s3) yield
{\small
\begin{align}
  \label{eq:tt}
  &\young(111,2221,33,4432)
  &
  &\young(111,:2221,:::33,::::4432)
  &
  &\young(:::::111,:::1222,:::33,2344)
  \\
  \tag*{}
  &\hphantom{\yng(2)}\Tt &&\hphantom{\yng(4)}\tshi
  &&\hphantom{\yng(2)} t
\end{align}
}%
We denote with $\Tt_{ij}$ (with $j=1, \ldots, \d_i$) the integer contained
in the $j^{\rm th}$ box of the $i^{\rm th}$ row of $\Tt$
(as if $\Tt$ were a matrix).
Recall that all we have to check is that the columns of $t$
are strictly increasing, since the other properties of the LR
tableaux (nondecreasing rows and the lattice permutation property)
are inherited from $\ho$. We are going to prove that

\begin{lemma}\label{th:ab}
\begin{enumerate}[(A)]
\item The columns of $\Tt$ are strictly increasing.%
\footnote{The columns of $\Tt$ in general are not contiguous, as in \eqref{eq:tt}, 
but it does not matter.}
\item Steps (s2) and (obviously) (s3) 
described above preserve the property
of having strictly increasing columns.
\end{enumerate}
\end{lemma}

\smallno
Proof of (A).
If $q$ is a positive integer, let
\begin{align*}
  h_\ell(q) &:= \text{number of 
  $q$'s appearing in the first $\ell$
  rows of $\Tt$.}
\end{align*}
We claim that
\begin{equation}\begin{split}
  \label{eq:pos}
  &\text{\itshape if $q$ appears in the $\ell^{\text{th}}$ row of $\Tt$,
  its first appearance} 
  \\[-1mm]
  &\text{\itshape 
  occurs at position $h_{\ell-1}(q)+1$, thus $q$
  occupies}
  \\[-1mm]
  &\text{\itshape
  the positions: 
  $\;h_{\ell-1}(q)+1,\: h_{\ell-1}(q)+2,\: \ldots,\:  h_{\ell}(q)$.}  
\end{split}\end{equation}
The second part of the statement trivially follows from the first
part and from the fact that $h_\ell(q) - h_{\ell-1}(q)$ 
is the number of $q$'s in the $\ell^{\rm th}$ row.
Let $p_1 > p_2 > \cdots > p_s$ be the distinct integers which appear
in the $\ell^{\text{th}}$ row of $\Tt$.
Thanks to Proposition \ref{th:ho}, we already know that $p_1 = \ell$
and that $h_{\ell-1}(\ell) = 0$, hence \eqref{eq:pos} holds for $p_1$.
Assume now that \eqref{eq:pos} holds for $p_i$. If we let 
$j = h_{\ell}(p_i)$, then $j$ is the last position where
$p_i$ appears, so
\begin{equation*}
  \ho_{z_\ell+j-1} = \Tt_{\ell, j} = p_i >  p_{i+1} = 
  \Tt_{\ell, j+1}  = \ho_{z_\ell+j}\,.
\end{equation*}
Statement (3) of Proposition \ref{th:ho} implies that 
$h_\ell(p_i) = h_{\ell-1}(p_{i+1})$.
By consequence,  $p_{i+1}$ makes its first appearance
in the $\ell^{\rm th}$ row of $\Tt$ at position 
\begin{equation*}
  j+1 = h_{\ell-1}(p_{i+1})+1 \,,
\end{equation*}
\iie\  \eqref{eq:pos} holds for $p_{i+1}$. By iteration we have that
\eqref{eq:pos} holds for all integers $p_m$ which are present
in the $\ell^{\rm th}$ row of $\Tt$.

\smallno
We can prove now that the columns of $\Tt$ are strictly increasing.
More precisely we prove the implication
\begin{equation}
  \label{eq:jj}
  \ell<m \quad\text{and}\quad\Tt_{\ell,j} \ge \Tt_{m,j'} 
  \quad\impp\quad j<j' \,.
\end{equation}
Let, in fact, $s\le \ol s$, and
\begin{equation*}
  \Tt_{\ell,j} = \ol s
  \qquad\text{and}\qquad
  \Tt_{m,j'} =   s \,.
\end{equation*}
From \eqref{eq:pos}, it follows that
\begin{equation*}
  j \le h_{\ell}(\ol s)
  \qquad\text{and}\qquad
  j' \ge h_{m-1}(s) + 1 \,.
\end{equation*}
Since $\ho$ is a lattice permutations, we get
$h_\ell(\ol s) \le h_\ell(s) \le h_{m-1}(s)$,
and \eqref{eq:jj} is proven.

\medno
Proof of (B). 
Using the definition $\tshi_{i,j} := \Tt_{i, j-\a_1+\a_i}$ and 
\eqref{eq:jj}, we obtain that, for $\ell<m$,
\begin{align*}
  \tshi_{\ell,j} \ge \tshi_{m,j'} 
  \quad\impp\quad
  j-\a_1+\a_\ell < j'-\a_1+\a_m
  \quad\impp\quad
  j<j'\,.
\end{align*}
Thus implication \eqref{eq:jj} holds for $\tshi$.
By consequence
the columns of $\tshi$ (and obviously those of $t$ as well)
are strictly increasing.
This concludes the proof of Lemma \ref{th:ab} and 
Proposition \ref{th:good}. Hence Lemma \ref{th:lrmin}
is also proven.
\qed

\secskip

%
%
%

\providecommand{\bysame}{\leavevmode\hbox to3em{\hrulefill}\thinspace}
\providecommand{\MR}{\relax\ifhmode\unskip\space\fi MR }
\providecommand{\MRhref}[2]{%
  \href{http://www.ams.org/mathscinet-getitem?mr=#1}{#2}
}
\providecommand{\href}[2]{#2}

\end{document}